\documentclass[preprint,12pt]{elsarticle}

\usepackage{graphicx}
\usepackage{amsmath} 
\usepackage{amsfonts} 
\usepackage{amssymb}
\usepackage{amsthm}  
\usepackage{latexsym}
\usepackage{eufrak}
\usepackage{euscript}
\usepackage{exscale}
\usepackage{graphicx}
\usepackage{mathrsfs}
\usepackage{hyperref}
\usepackage{stmaryrd}
\usepackage{xypic}
\usepackage{color}
\usepackage[margin=2.5cm]{geometry}
\newcommand{\kc}{\mathcal{C}}

\newcommand{\pc}{\mathcal{P}_\mathcal{C}}

\journal{}
\makeatletter
\def\ps@pprintTitle{%
 \let\@oddhead\@empty
 \let\@evenhead\@empty
 \def\@oddfoot{}%
 \let\@evenfoot\@oddfoot}

\begin{document}

\begin{frontmatter}

\title{On distributions with fixed marginals maximizing the joint or the prior default probability, estimation,
 and related results}

\author[add1]{Thomas Mroz}
\ead{thomasmroz@a1.net}
\author[add2]{Juan Fern\'andez S\'anchez}
\ead{juanfernandez@ual.es}
\author[add1]{Sebastian Fuchs}
\ead{sebastian.fuchs@plus.ac.at}
\author[add1]{Wolfgang Trutschnig\corref{cor1}}
\ead{wolfgang@trutschnig.net}
\cortext[cor1]{Corresponding author}

\address[add1]{Department of Mathematics, University of Salzburg, Hellbrunner Strasse 34, 5020 Salzburg, Austria,\\
Tel.: +43 662 8044-5312,  Fax: +43 662 8044-137}
\address[add2]{Grupo de Teor\'ia de C\'opulas y aplicaciones, Universidad de Almer\'ia, La Ca\~nada de San Urbano, Almer\'ia, Spain }   

\begin{abstract}
We study the problem of maximizing the probability that (i) an electric component or financial institution $X$ does not 
default before another component or institution $Y$ and (ii) that $X$ and $Y$ default jointly within the class of all random variables $X,Y$ with given univariate continuous 
distribution functions $F$ and $G$, respectively, and show that the maximization problems correspond to 
finding copulas maximizing the mass of the endograph $\Gamma^\leq(T)$ and the graph $\Gamma(T)$ of 
$T=G \circ F^-$, respectively. 
After providing simple, copula-based proofs for the existence of copulas attaining the two maxima $\overline{m}_T$ and
$\overline{w}_T$ we generalize the obtained results to the case of general (not necessarily monotonic) transformations 
$T:[0,1] \rightarrow [0,1]$ and derive simple and easily calculable formulas for $\overline{m}_T$ and $\overline{w}_T$
involving the distribution function $F_T$ of $T$ (interpreted as random variable on $[0,1]$).
The latter are then used to charac\-terize all non-decreasing transformations $T:[0,1] \rightarrow [0,1]$ for which
$\overline{m}_T$ and $\overline{w}_T$ coincide. A strongly consistent estimator for the maximum probability that
$X$ does not default before $Y$ is derived and proven to be asymptotically normal under very mild regularity conditions.
Several examples and graphics illustrate the main results and falsify some seemingly natural conjectures.
\end{abstract}

\begin{keyword}
Copula \sep Dependence \sep Estimator \sep Graph \sep Endograph \sep Markov Kernel\\[0.4cm]
 \MSC[2010] 60E05 \sep 28A50 \sep 91G70 
\end{keyword}

\end{frontmatter}

 
\setlength{\textwidth}{480pt} 

\newtheorem{thm}{Theorem}
\newtheorem{lem}[thm]{Lemma}
\newtheorem{prop}[thm]{Proposition}
\newtheorem{cor}[thm]{Corollary}
\newdefinition{rmk}[thm]{Remark}
\newdefinition{ex}[thm]{Example}
\newdefinition{defin}[thm]{Definition}

\section{Introduction}
Suppose that $F$ and $G$ are (continuous) distribution functions of two random variables $X$ and $Y$ modeling, e.g., 
(i) the default times of financial institutions or (ii) the lifetime of electronic components. 
Especially in the context of (i)
the marginal distributions might be know or at least be estimated in standard ways, whereas 
the joint distribution is often unknown and harder to estimate. In such situations (particularly in the context of so-called credit default swaps) is seems natural to consider the worst-case scenario and study bivariate distribution 
functions $H$ in the Fr\'echet class $\mathcal{H}_{F,G}$ of $F,G$  (the family of all bivariate distribution functions 
having marginals $F$ and $G$) with the following property: In case $(X,Y)$ has distribution function $H$ 
the joint or prior default probability (i.e. the probability of the events $\{X=Y\}$ and $\{X \geq Y\}$, respectively) is maximal within $\mathcal{H}_{F,G}$.

Translating to the class of copulas (see \cite{MS} and Section 2), maximizing the afore-mentioned probabilities 
means calculating 
\begin{equation}\label{eq:objective}
\overline{w}_T:=\sup_{A \in \kc} \mu_A(\Gamma(T)), \quad \overline{m}_T:=\sup_{A \in \kc} \mu_A(\Gamma^\leq(T))
\end{equation} 
where $T:[0,1] \rightarrow [0,1]$ is defined by $T=G \circ F^- $, $F^-$ denotes the quasi-inverse of $F$, $\Gamma(T)$ the graph of $T$, $\Gamma^\leq(T)=\{(x,y)\in [0,1]^2: y \leq T(x)\}$ the so-called endograph of $T$, $\kc$ the family of all two-dimensional copulas and $\mu_A$ the doubly stochastic measure corresponding to the copula $A \in \kc$.  

It has been brought to our attention that formulas for the suprema in eq. (\ref{eq:objective}) also follow 
from deep and much heavier machinery going back to Rüschendorf in \cite{Rue}. 
In the current paper we provide (a) independent alternative simple, copula-based proofs and show the existence of 
copulas $A_R, A_S \in \kc$ attaining the suprema (including the fact that
it is possible to choose $A_R$ completely dependent). Complementing these results, (b) we calculate 
$\overline{w}_T$ and $\overline{m}_T$ also for general measurable, not necessarily monotonic transformations
$T:[0,1] \rightarrow [0,1]$, (c) characterize for which non-decreasing $T$ we even have
$\overline{w}_T=\overline{m}_T$ and, (d) derive a strongly consistent estimator for 
$\overline{m}_T$ and show that the latter is asymptotically normal under mild regularity conditions.    
 
The rest of the paper is organized as follows: Section \ref{sec:notation} gathers some preliminaries and notations, and
proves the afore-mentioned translation of the problem of maximizing the joint or prior default probability 
to the copula setting. 
The main results concerning the calculation of the maximum probabilities and various related questions 
are gathered in Sections \ref{sec:maxendograph} and \ref{sec:maxgraph}, whereas in Section \ref{sec:coincide} we characterize the case $\overline{w}_T=\overline{m}_T$ for
non-decreasing $T$. Finally, Section \ref{sec:asymptotics} introduces an estimator for $\overline{m}_T$, shows 
consistency and studies its asymptotic distribution. Several examples and graphics illustrate the obtained results and the chosen approach. 

\section{Notation and Preliminaries}\label{sec:notation}
For every $d$-dimensional random vector $\mathbf{X}$ on a probability space $(\Omega,\mathcal{A},\mathbb{P})$ we will write $\mathbf{X} \sim F$ 
if $\mathbf{X}$ has distribution function (d.f., for short) $F$ and let $\mu_F=\mathbb{P}^{\mathbf{X}}$ denote the corresponding distribution on the 
Borel $\sigma$-field $\mathcal{B}(\mathbb{R}^d)$ of $\mathbb{R}^d$.
For every univariate distribution function $F$ we will let $F^-$ denote the quasi-inverse of $F$, i.e.
$F^-(q)=\inf\{x \in \mathbb{R}: F(x) \geq q\}$. Note that for every $q \in (0,1)$ we have 
$F^-(q) \leq x$ if and only if $q \leq F(x)$, that for $X \sim F$ and $F$ continuous we have $F \circ X \sim \mathcal{U}(0,1)$, i.e., $F \circ X$  is uniformly distributed on $[0,1]$, and that 
the random variable $F^-\circ F \circ X$ coincides with $X$ with probability one. For further properties of $F^-$ we refer, for instance, to \cite{EH}.
Given univariate distribution functions $F$ and $G$, we will let $\mathcal{H}_{F,G}$ denote the Fr\'echet class of $F$ and $G$, i.e. the family of all 
two-dimensional distribution functions having $F$ and $G$ as marginals; $\mathcal{P}_{F,G}$ will denote the corresponding class of probability measures on 
$\mathcal{B}(\mathbb{R}^2)$. $\mathcal{B}([0,1])$ and $\mathcal{B}([0,1]^2)$ denote the Borel $\sigma$-fields on $[0,1]$ and $[0,1]^2$, 
$\lambda$ and $\lambda_2$ the Lebesgue measure on $\mathcal{B}([0,1])$ and $\mathcal{B}([0,1]^2)$ respectively. For every measurable transformation $T:[0,1] \rightarrow [0,1]$
the push-forward of $\lambda$ via $T$ will be denoted by $\lambda^T$, i.e., $\lambda^T(E)=\lambda(T^{-1}(E))$ for every
$E \in \mathcal{B}([0,1])$. 

As already mentioned before, $\kc$ will denote the family of all two-dimen\-sional \emph{copulas}. For background on copulas we refer to \cite{Du,Nel}.
$M$ and $W$ will denote the upper and lower Fr\'echet-Hoeffding bounds, $\Pi$ the product copula.  
$d_\infty$ will denote the uniform distance on $\kc$; it is well known that $(\kc,d_\infty)$ is a compact metric space and that
$d_\infty$ is a metrization of weak convergence in $\kc$. 
For every $A \in \kc$ 
$\mu_A$ will denote the corresponding \emph{doubly stochastic measure} defined via 
$\mu_A([0,x]\times [0,y])=A(x,y)$ for all $x,y\in [0,1]$ (and extended in the standard way to $\mathcal{B}([0,1]^2)$), 
$\mathcal{P}_\mathcal{C}$ the class of all these doubly stochastic measures.

A \emph{Markov kernel} from $\mathbb{R}$ to $\mathcal{B}(\mathbb{R})$ is a mapping 
$K: \mathbb{R} \times \mathcal{B}(\mathbb{R})\rightarrow [0,1]$ 
such that $x \mapsto K(x,B)$ is measurable for every fixed $B \in \mathcal{B}(\mathbb{R})$ and $B \mapsto K(x,B)$ is a 
probability measure for every fixed $x \in \mathbb{R}$. 
Given real-valued random variables $X,Y$ on $(\Omega, \mathcal{A}, \mathbb{P})$, a Markov kernel 
$K:\mathbb{R}\times \mathcal{B}(\mathbb{R}) \rightarrow [0,1]$ is called a \emph{regular conditional distribution of 
$Y$ given $X$} if for every $B \in \mathcal{B}(\mathbb{R})$ 
\begin{equation}
  K(X(\omega),B)=\mathbb{E}(\mathbf{1}_B\circ Y |X)(\omega) 
\end{equation} 
holds $\mathbb{P}$-a.s. It is well known that for each pair $(X,Y)$ of real-valued random variables a regular 
conditional distribution $K(\cdot,\cdot)$ of $Y$ given $X$ exists, that $K(\cdot,\cdot)$ is unique $\mathbb{P}^X$-a.s.
(i.e. unique for $\mathbb{P}^X$-almost every $x \in \mathbb{R}$) and that $K(\cdot,\cdot)$ only depends on the distribution $\mathbb{P}^{(X,Y)}$. 
Hence, given $(X,Y)\sim H$, we will denote (a version of) the regular conditional distribution of $Y$ given $X$ 
by $K_H(\cdot,\cdot)$ and refer to $K_H(\cdot,\cdot)$ simply as \emph{Markov kernel of $H$} or \emph{Markov kernel of $(X,Y)$}.
Note that for every two-dimensional distribution function $H$, its Markov kernel $K_H(\cdot,\cdot)$, and every Borel set 
$G \in  \mathcal{B}(\mathbb{R}^2)$ the following \emph{disintegration} formula holds ($G_x=\{y \in \mathbb{R}:(x,y) \in G\}$ denoting the $x$-section of $G$ for 
every $x \in \mathbb{R}$)
\begin{equation}\label{disallg}
 \int_{\mathbb{R}} K_H(x,G_x)\, d\lambda(x) = \mu_H(G).
\end{equation}
For $A \in \kc$ we will directly consider the corresponding Markov kernel $K_A(\cdot,\cdot)$ to be defined on $[0,1] \times \mathcal{B}([0,1])$. Considering that
in this case eq. (\ref{disallg}) implies that
\begin{equation}\label{bound2}
 \int_{[0,1]} K_A(x,F)\, d\lambda(x) = \lambda(F)
\end{equation}
holds for every $F \in \mathcal{B}([0,1])$, and that, additionally, every Markov kernel $K:[0,1]\times \mathcal{B}([0,1]) \rightarrow [0,1]$ fulfilling eq. (\ref{bound2}) 
obviously induces a unique element $\mu \in \pc$, it follows that there is a one-to-one correspondence between $\kc$ and the family of 
all Markov kernels $K:[0,1] \times \mathcal{B}([0,1]) \rightarrow [0,1]$ fulfilling eq. (\ref{bound2}). Notice that for $A\in\kc$ eq. (\ref{bound2}) also implies that
$K_A(x,\{0,1\})=0$ holds for $\lambda$-almost every $x \in [0,1]$, so it is always possible to choose a (version of the) kernel fulfilling $K_A(x,\{0,1\})=0$ for every $x \in [0,1]$.
For more details and properties of conditional expectation, regular conditional distributions, and disintegration 
see \cite{Ka} and \cite{Kl}, various results underlining the usefulness of the Markov kernel perspective can be found in \cite{Du} and the references therein.

In the sequel $\mathcal{T}$ will denote the class of all $\lambda$-preserving transformations $h:[0,1] \rightarrow [0,1]$, i.e., the class of all $h$ fulfilling $\lambda^h=\lambda$,  
$\mathcal{T}_b$ the subset of all bijective $h \in \mathcal{T}$, and $\mathcal{T}_l$ the subset of all piecewise linear, bijective $h \in \mathcal{T}$. 
A co\-pula $A\in \kc$ will be called \emph{completely dependent} if and only if there exists $h \in \mathcal{T}$ such that $K(x,E)=\mathbf{1}_E(h(x))$ is a 
regular conditional distribution of $A$ (see \cite{Lan,Tru} for equivalent definitions and main properties). 
For every $h \in \mathcal{T}$ the induced completely dependent copula will be denoted by $A_h$ throughout the rest of the 
paper, $\kc_d$ will denote the family of all completely dependent copulas.

Following \cite{Du,TFS}, for every $h \in \mathcal{T}$ and every copula $A\in\kc$ we will let
$\mathcal{S}_h(A) \in \kc$ denote the (generalized) \emph{$h$-shuffle of $A$}, defined implicitly via the corresponding doubly stochastic measures by
\begin{equation}\label{defshuffle}
\mu_{\mathcal{S}_h(A)}(E \times F)=\mu_A(h^{-1}(E) \times F)
\end{equation}
for all $E,F \in \mathcal{B}([0,1])$. Notice that $\mathcal{S}_h(A)$ is a shuffle in the sense of \cite{DSS} if $h \in \mathcal{T}_b$, and that for $A=M$ it is a shuffle in the 
sense of \cite{MST} (to which we will refer as classical shuffle in the sequel) if $h \in \mathcal{T}_l$.

We conclude this section with the afore-mentioned translation of the maximization problems to the copula setting and start
with the following lemma which is straightforward to prove via disintegration and a Dynkin system argument. 

\begin{lem}\label{lem:KAtoKH}
Suppose that $F,G$ are continuous distribution functions, that $(X,Y)$ has d.f. $ H \in \mathcal{H}_{F,G}$ and copula $A$, and let 
$K_A(\cdot,\cdot)$ denote a Markov kernel of $A$ fulfilling $K_A(x,\{0,1\})=0$ for all $x \in [0,1]$. Then setting
\begin{equation}\label{eq:defKH}
K\big(x,(-\infty,y]\big) := K_A\big(F(x),[0,G(y)] \big)
\end{equation}
for all $x,y \in \mathbb{R}$ defines a Markov kernel $K(\cdot,\cdot)$ of $(X,Y) \sim H$.
\end{lem}
Suppose now that $S:\mathbb{R} \rightarrow \mathbb{R}$ is an arbitrary Borel-measurable mapping. In the sequel we will let $\Gamma(S)$ and $\Gamma^{\leq}(S)$ denote the graph
and the endograph of $S$ respectively, i.e.
\begin{equation}
\Gamma(S)=\{(x,S(x)): x \in \mathbb{R}\}, \quad \Gamma^{\leq}(S)=\{(x,y) \in \mathbb{R}^2: y \leq S(x)\}.
\end{equation}
Lemma \ref{lem:KAtoKH} allows to express $\mathbb{P}(Y\leq X)$ as well as $\mathbb{P}(Y=X)$ in terms of $F,G$ and the underlying copula $A$. 
In order to prove a more general result and to simplify notation, given (continuous) $F,G$ and (measurable) $S$ we will write 
\begin{equation}\label{ST}
T:=G \circ S \circ F^-
\end{equation} 
in the sequel. In general, $T$ is only well-defined on $(0,1)$ - we will however, directly consider it as function on $[0,1]$ by setting $T(0):=0$ and $T(1):=T(1-)= \lim_{x \rightarrow 1-}T(x)$. 
\begin{thm}\label{translate}
Suppose that $X,Y$ are random variables on $(\Omega,\mathcal{A},\mathbb{P})$ with joint distribution function $H$, continuous marginals $F$ and $G$ and copula $A$. 
Furthermore let $S:\mathbb{R} \rightarrow \mathbb{R}$ be an arbitrary Borel-measurable mapping and define $T$ according to eq. (\ref{ST}). 
Then the following identities hold for $T:=G \circ S \circ F^-$:
\begin{equation}
\mathbb{P}^{(X,Y)}\big(\Gamma(S)\big) = \mu_A(\Gamma(T)),\quad \mathbb{P}^{(X,Y)}\big(\Gamma^{\leq}(S)\big) = \mu_A(\Gamma^{\leq}(T))
\end{equation}
\end{thm}
\begin{proof} Using the fact that $\mathbb{P}(F^-\circ F \circ X=X)=1$, change of coordinates, 
disintegration and Lemma \ref{lem:KAtoKH} the second identity can be proved as follows:
\begin{eqnarray*}
\mathbb{P}^{(X,Y)}\big(\Gamma^{\leq}(S)\big) &=& \int_\Omega K_H \big(X(\omega),(-\infty,S \circ X (\omega)] \big) \,d\mathbb{P}(\omega) \\
    &=& \int_\Omega K_A \big(F \circ X(\omega),[0,G \circ S \circ F^- \circ F \circ X (\omega)] \big) \,d\mathbb{P}(\omega) \\
    &=& \int_{[0,1]} K_A \big(z,[0,G \circ S \circ F^- (z)] \big) \,d\mathbb{P}^{F \circ X}(z) \\
    &=&  \int_{[0,1]} K_A \big(z,[0,T(z)] \big) \,d\lambda(z) = \mu_A(\Gamma^{\leq}(T)).
\end{eqnarray*}
Working with $K \big(X(\omega),\{S \circ X (\omega)\} \big)$ instead of $K \big(X(\omega),(-\infty,S \circ X (\omega)] \big)$ the first identity
$\mathbb{P}^{(X,Y)}\big(\Gamma(S)\big) = \mu_A(\Gamma(T))$ follows in the same manner. 
\end{proof}

\section{Maximizing the mass of the endograph and the prior default probability} \label{sec:maxendograph}
\begin{figure}[h!]
  \begin{center}
  \includegraphics[width = 10cm,angle=270]{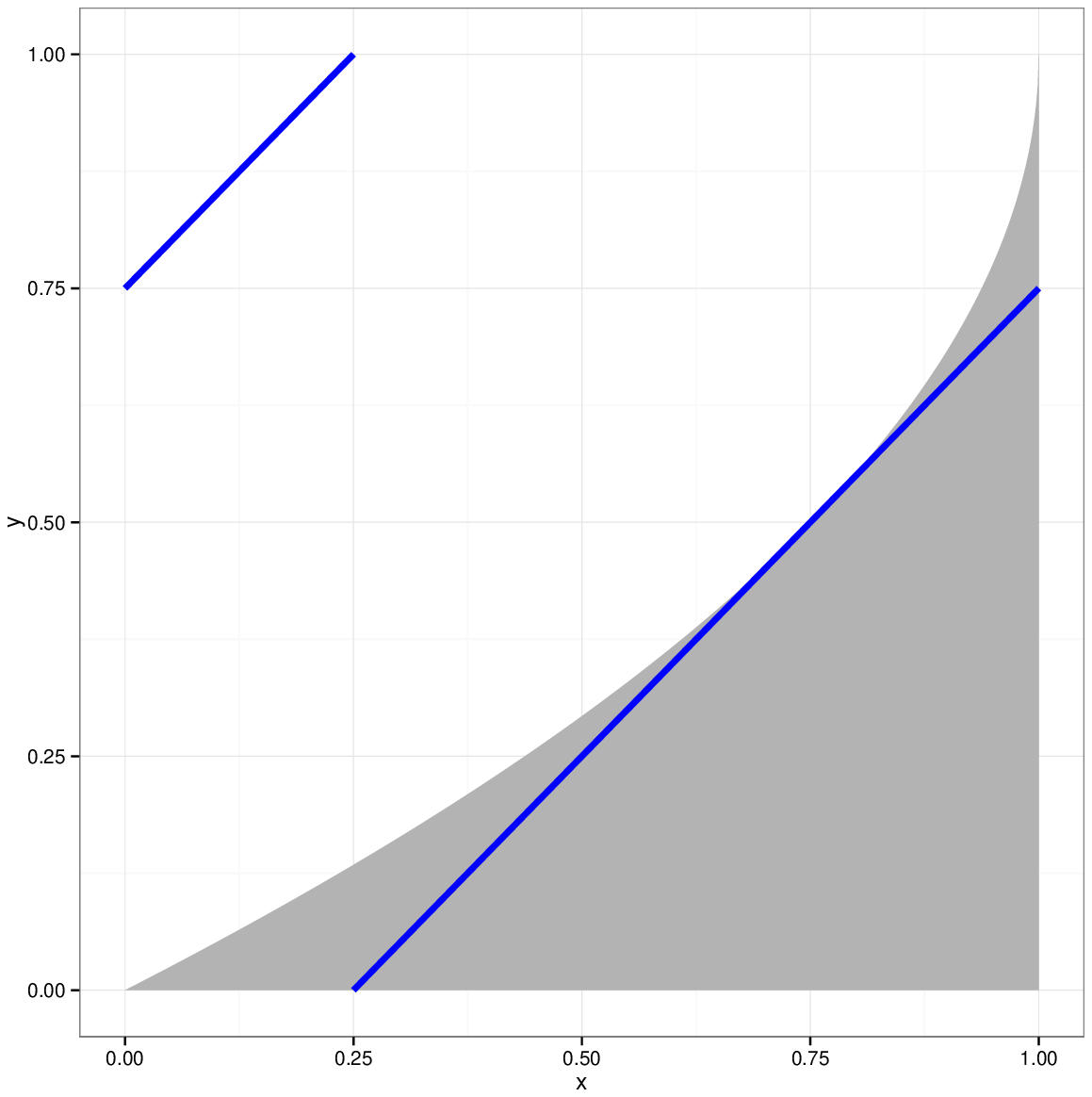}
        \caption{The endograph $\Gamma^\leq(T)$ of the transformation $T(x)=1-(1-x)^{\frac{1}{2}}$ (shaded region) and 
        the support of the mutually completely dependent copula $A_R$ constructed 
        in the proof of Theorem \ref{solutionmax} assigning maximum mass to $\Gamma^\leq(T)$ (blue).}\label{grx2}
  \end{center}
\end{figure}

Suppose that $X \sim F$ and $Y \sim G$ model default times and that $F,G$ are continuous. 
Considering $S=id_\mathbb{R}$ then calculating 
$\sup_{\mu \in \mathcal{P}_{F,G}} \mu(\Gamma^{\leq}(S))$ obviously corresponds to finding (joint) distributions of 
$(X,Y)$ maximizing the probability of a prior or joint default. To simplify notation in the sequel we will 
simply refer to the event $\{Y \leq X\}$ as `prior default' (of $Y$) although $\{Y \leq X\}$ corresponds to 
the prior and joint default. 
Notice that, setting $\psi(x,y)=x+y$ and considering the pair $(-X,Y)$ the afore-mentioned maximization problem can 
be considered a special case of the more general situation studied in \cite{EP,EP2}.
Theorem \ref{translate} implies
\begin{equation}\label{maxtranslate}
\overline{m}_{F,G}:=\sup_{\mu \in \mathcal{P}(F,G)} \mu(\Gamma^{\leq}(id_\mathbb{R})) = \sup_{A \in \kc} \mu_A (\Gamma^{\leq}(T))=:\overline{m}_T
\end{equation}  
as well as
\begin{equation}\label{mintranslate}
\underline{m}_{F,G}:=\inf_{\mu \in \mathcal{P}(F,G)} \mu(\Gamma^{\leq}(id_\mathbb{R})) = \inf_{A \in \kc} \mu_A (\Gamma^{\leq}(T))=:\underline{m}_T
\end{equation} 
whereby $T=G \circ S \circ F^- = G \circ F^-$. Since $G \circ F^-$ is non-decreasing it is possible to 
derive a simple formula for $\underline{m}_T$ and even construct a dependence structure for which $\mathbb{P}(Y \leq X)$
coincides with $\overline{m}_T$. The following result holds:
\begin{thm} \label{solutionmax}
Suppose that $T:[0,1] \rightarrow [0,1]$ is non-decreasing. Then we have 
\begin{equation}\label{maxendo}
\overline{m}_T= \sup_{A \in \kc} \mu_A\big(\Gamma^{\leq}(T) \big) = 1+\inf_{x \in [0,1]} (T(x)-x).
\end{equation} 
Moreover, defining $R \in \mathcal{T}$ by $R(x)=x+\overline{m}_T \,(mod\, 1)$, 
we have \mbox{$\mu_{A_R}(\Gamma^{\leq}(T))=\overline{m}_T$.}
\end{thm}
\begin{proof} Considering $\Gamma^{\leq}(T) \subseteq [0,x] \times [0,T(x)] \, \cup \, [x,1] \times [0,1]$ it follows that 
$\mu_A(\Gamma^{\leq}(T)) \leq T(x) + 1-x$ holds for every $x \in [0,1]$ and every $A \in \kc$, which implies that the left-hand side of (\ref{maxendo}) 
is smaller than or equal to the right-hand side. \\
To prove the reverse inequality set $z=\inf_{x \in [0,1]} \big(T(x)+1-x \big)$. 
For $z=1$ we have $T(x) \geq x$ for every $x$, so  taking into account $\mu_M(\Gamma^{\leq}(T))=1$ we are done, and 
it suffices to consider $z <1$. Compactness of $[0,1]$ implies the existence of a sequence 
$(x_n)_{n \in \mathbb{N}}$ and a point $x^\star \in [0,1]$ such that $\lim_{n \rightarrow \infty} x_n=x^\star$ and $\lim_{n \rightarrow \infty} (T(x_n)+1-x_n)=z$.
Using $z<1$ we get $x^\star >0$ and, using monotonicity of $T$ 
it follows that $T(x^\star-) + 1-x^\star=z$. Letting $R:[0,1] \rightarrow [0,1]$ denote the rotation defined by 
$R(x)=x+z \,(mod\, 1)$, obviously $R \in \mathcal{T}$ holds. Considering that for every $x \in [x^\star - T(x^\star-), 1]$ we have (see Figure \ref{grx2})
\begin{eqnarray*}
R(x)&=& T(x^\star-)-x^\star + x = T(x^\star-)+1-x^\star -1 + x \\
    &\leq & T(x)+1-x -1 +x =T(x)
\end{eqnarray*} 
it follows immediately that 
\begin{eqnarray*}
\mu_{A_R}(\Gamma^{\leq}(T)) &\geq& 1 - (x^\star - T(x^\star-)) =z=\inf_{x \in [0,1]} \big(T(x)+1-x \big), 
\end{eqnarray*}
which completes the proof. 
\end{proof}
\begin{rmk}
Considering that continuity of $T$ plays no role in Theorem \ref{solutionmax}, that $T$ has (as non-decreasing function) at most countably many
discontinuities, and that $\mu_A(E \times [0,1])=0$ for every countable set $E$ and $A \in \kc$ we may, w.l.o.g., assume that $T$ is left continuous, in which case the 
infimum in eq. (\ref{maxendo}) is a minimum.
\end{rmk}
\begin{cor}
Suppose that $X,Y$ are random variables with continuous distribution functions $F$ and $G$ respectively, set $T=G \circ F^-$ and 
$z:=1+\inf_{x \in [0,1]} (T(x)-x)$, define $R:[0,1] \rightarrow [0,1]$ by $R(x)=z+x \,(mod\, 1)$, and let $A_R$ denote the completely dependent copula induced by $R$. 
Then for $(X,Y) \sim H \in \mathcal{H}(F,G)$ with $H(x,y)=A_R(F(x),G(y))$ we have $\mathbb{P}(Y \leq X)= 
\overline{m}_{F,G}$.
\end{cor}
\begin{ex}\label{exponential_lifetimes}
Suppose that the default times $X$ and $Y$ are exponentially distributed with parameters $\theta_1$ and $\theta_2$,  respectively. It is straightforward to verify that in this case
$T=G \circ F^-$ is given by $T_\theta(x)=1-(1-x)^\theta$, where $\theta=\frac{\theta_2}{\theta_1}$. For the case of $\theta \geq 1$ we have $T_\theta(x)\geq x$ for every 
$x \in [0,1]$, so $ \sup_{A\in \kc} \mu_A(\Gamma^\leq(T_\theta))=1$. 
Remarkably, for the case of $\theta<1$ the maximal mass of the endograph of $T_\theta$ and the maximal mass of the graph of $T_\theta$ coincide. In fact, 
applying Theorem \ref{solutionmax}, on the one hand we get 
$$
\sup_{A \in \kc} \mu_A\big(\Gamma^{\leq}(T_\theta) \big)= 1 + \theta^{\frac{1}{1-\theta}} - \theta^{\frac{\theta}{1-\theta}}.
$$
And on the other hand, according to Theorem 3 and Theorem 4 in \cite{DFST} (also see \cite{MS, Th}) we have 
\begin{equation}\label{maxgraphjap_exponential}
\sup_{A\in \kc} \mu_A(\Gamma(T_\theta))= \int_{[0,1]} \Big(\mathbf{1}_{[0,1]} (f \circ T_\theta) + \frac{1}{f\circ T_\theta} \mathbf{1}_{(1,\infty)}(f\circ T_\theta)\Big) \,d\lambda
\end{equation}
where $f$ denotes the density of $\lambda^{T_\theta}$. Since the latter is given by 
$f(x)=\frac{1}{\theta}\,(1-x)^{\frac{1-\theta}{\theta}}$ we get $f\circ T_\theta(x)= \frac{1}{\theta} (1-x)^{1-\theta}$ and 
eq. (\ref{maxgraphjap_exponential}) calculates to
\begin{eqnarray*}
\sup_{A\in \kc} \mu_A(\Gamma(T_\theta)) &=& \int_{\big[0,1-\theta^{\frac{1}{1-\theta}}\big]} \frac{1}{\frac{1}{\theta} (1-x)^{1-\theta}} d\lambda(x) + 
     1-\big(1- \theta^{\frac{1}{1-\theta}}\big)  
 =  1 - \theta^{\frac{\theta}{1-\theta}} + \theta^{\frac{1}{1-\theta}}\\
  &=& \sup_{A \in \kc} \mu_A\big(\Gamma^{\leq}(T_\theta) \big).
\end{eqnarray*}
For the special case of $\theta=\frac{1}{2}$, which is depicted in Figure \ref{grx2}, we get 
$$
\sup_{A\in \kc} \mu_A(\Gamma(T))=\sup_{A \in \kc} \mu_A\big(\Gamma^{\leq}(T) \big)=\frac{3}{4}.
$$
\end{ex}

\begin{figure}[h!]
  \begin{center}
  \includegraphics[width = 10cm,angle=270]{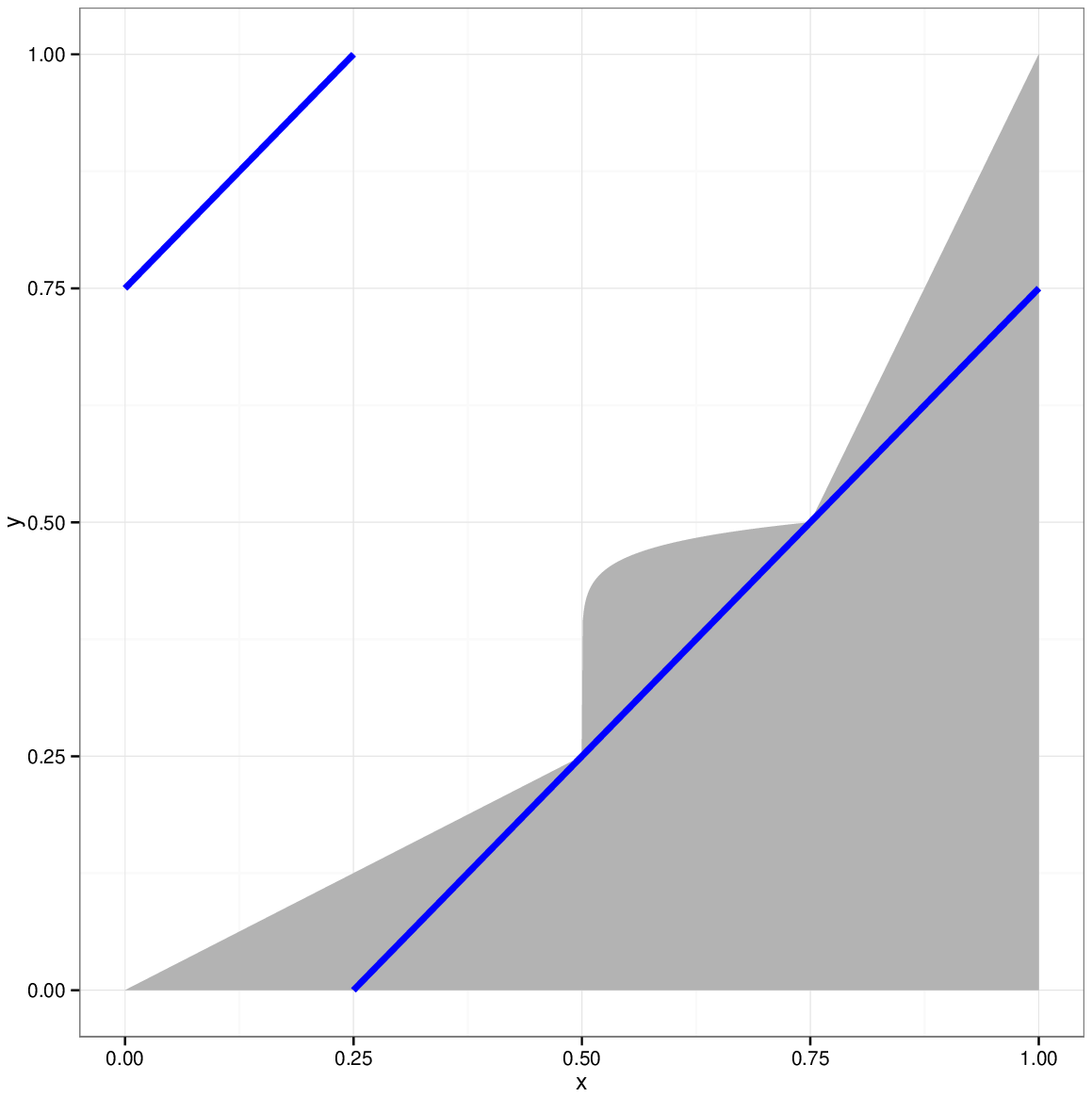}
        \caption{The endograph $\Gamma^\leq(T_{10})$ of the transformation $T_{10}$ from Example \ref{exgegen} (shaded region) and the 
        support of the mutually completely dependent copula $A_R$ constructed 
        in the proof of Theorem \ref{solutionmax} assigning maximum mass to $\Gamma^\leq(T_{10})$ (blue).}\label{homgegenbsp}
  \end{center}
\end{figure}

\begin{ex}\label{exgegen}
Based on Example \ref{exponential_lifetimes} it might seem natural to conjecture that the equality 
$\sup_{A\in \kc} \mu_A(\Gamma(T))=\sup_{A \in \kc} \mu_A\big(\Gamma^{\leq}(T) \big)$ holds
for a much bigger class of non-decreasing transformations $T$ fulfilling $T(x)\leq x$ for every $x \in [0,1]$. 
Since counterexamples are easily constructed for the case where $T$ is singular ($\lambda^T(E)>0$ for some $E \in \mathcal{B}([0,1])$ with $\lambda(E)=0$) 
and the case where $T$ has discontinuities, the conjecture reduces to strictly increasing, continuous transformations $T$. For every $n \in \mathbb{N}$ 
the transformation $T_n:[0,1] \rightarrow [0,1]$, defined by 
$$
T_n(x) = \left\{
\begin{array}{rl}
\frac{x}{2} & \textrm{if } x \in [0,\frac{1}{2}] \\
\frac{x}{2} + \frac{x}{2} \sqrt[n]{4x-2} & \textrm{if } x \in (\frac{1}{2},\frac{3}{4}) \\
-1+2x & \textrm{if } x \in [\frac{3}{4},1]
\end{array} \right.
$$ 
is easily verified to be homeomorphism with $T_n (x)\leq x$ for every $x \in [0,1]$ (see Figure \ref{homgegenbsp} for the case $n=10$).
Applying Theorem \ref{solutionmax} we get $\sup_{A\in \kc} \mu_A(\Gamma^\leq(T))=\frac{3}{4}$, however, either by graphical arguments or by 
using Theorem 3 and Theorem 4 in \cite{DFST}
it is straightforward to verify that $\lim_{n \rightarrow \infty}\sup_{A\in \kc} \mu_A(\Gamma(T_n))=\frac{1}{2} < \frac{3}{4}$, so the conjecture is wrong.
\end{ex}
Although monotonicity is crucial in the proof of Theorem \ref{solutionmax} it is even possible to calculate 
\begin{equation}
\overline{m}:=\sup_{\mu \in \mathcal{P}(F,G)} \mu(\Gamma^{\leq}(S)) = \sup_{A \in \kc} \mu_A(\Gamma^{\leq}(T)) \nonumber
\end{equation}
for the case of arbitrary measurable (not necessarily monotonic) transformations $S: \mathbb{R} \rightarrow \mathbb{R}$ 
(as before $T:=G \circ S \circ F^-$). 
Letting $T:[0,1] \rightarrow [0,1]$ denote an arbitrary measurable transformation, we will now directly concentrate on the quantity 
\begin{equation}\label{MT}
\overline{m}_T:= \sup_{A \in \kc} \mu_A(\Gamma^{\leq}(T))
\end{equation}
and prove a simple formula for $\overline{m}_T$ only involving the d.f. $F_T:[0,1] \rightarrow [0,1]$ of $T$, defined by
\begin{equation}
F_T(x)=\lambda^T([0,x])=\lambda(T^{-1}([0,x])).
\end{equation} 
We start with two simple lemmata that will be used in the proof of the main results.
\begin{lem}\label{lemupper}
Suppose that $T:[0,1] \rightarrow [0,1]$ is measurable. Then we have
\begin{equation}\label{maxendogen}
  \overline{m}_T 
	\leq 
	     1 + \inf_{y \in [0,1]} \big(y-F_T(y)\big) =1 + \min_{y \in [0,1]} \big(y-F_T(y)\big)
\end{equation} 
If $T$ is non-decreasing then we have equality in (\ref{maxendogen}). 
\end{lem}
\begin{proof}
Considering $\Gamma^{\leq}(T) \subseteq [0,1] \times [0,y] \, \cup \, T^{-1}((y,1]) \times [0,1]$ and using 
$\lambda^T((y,1])=1-F_T(y) $ we get   
$$
\mu_A(\Gamma^{\leq}(T)) \leq y + 1- F_T(y)
$$
for every $y \in [0,1]$ and every $A \in \kc$, from which the first inequality follows immediately.\\
Proving the existence of $y^\star \in [0,1]$ fulfilling $I:=\inf_{y \in [0,1]} (y-F_T(y)) =  y^\star-F_T(y^\star)$ can be done as follows: For every 
$n \in \mathbb{N}$ we can find $y_n \in [0,1]$ with  $y_n-F_T(y_n) < I + \frac{1}{2^n}$. Compactness of $[0,1]$ implies the existence of a subsequence $(y_{n_j})_{j \in \mathbb{N}}$
and some $y^\star \in [0,1]$ with $\lim_{j \rightarrow \infty} y_{n_j} =y^\star$. If $y^\star=1$ we are done since 
$I=\lim_{j \rightarrow \infty} (y_{n_j}-F_T(y_{n_j})) = y^\star- \lim_{j \rightarrow \infty} F_T(y_{n_j}) \geq y^\star - 1 = y^\star - F_T(y^\star)$. Suppose therefore that 
$y^\star < 1$ and let $\delta \in (0,1-y^\star]$ be arbitrary. Then there exists an index $j_0 \in \mathbb{N}$ such that $y_{n_j}<y^\star + \delta$, hence
$y_{n_j}-F_T(y_{n_j}) \geq y_{n_j} - F_T(y^\star + \delta)$, holds for all $j \geq j_0$. Considering $j \rightarrow \infty$ yields $I \geq y^\star - F_T(y^\star + \delta)$, hence,
using right-continuity of $F_T$ we get $I \geq y^\star - F_T(y^\star)$.   \\
Finally, suppose that $T$ is non-decreasing. We want to show that 
\begin{equation}\label{mono2}
\inf_{y \in [0,1]} (y-F_T (y)) = \inf_{x \in [0,1]} (T(x)-x)
\end{equation}
It follows directly from the construction that $F_T\circ T(x) \geq x$ holds for every $x \in [0,1]$
implying
$$
  \inf_{y \in [0,1]} (y-F_T (y)) 
	\leq T(x)- F_T(T(x))
	\leq T(x)- x
$$
for every $x \in [0,1]$ and hence
$$
  \overline{m}_T
	\leq 1 + \inf_{y \in [0,1]} \big(y-F_T(y)\big)
	\leq 1 + \inf_{x \in [0,1]} (T(x)-x)
	  =  \overline{m}_T
$$
which completes the proof. 
\end{proof}

\begin{lem}\label{sensitiv}
Suppose that $T,T':[0,1] \rightarrow [0,1]$ are measurable transformations. Then the following two assertions hold:
\begin{enumerate}
\item For $D:=\{x \in [0,1]:T(x)\not = T'(x)\}$ we have $\vert\overline{m}_{T'}- \overline{m}_T\vert \leq \lambda(D)$.
\item If $\Delta\in [0,1)$ and $T' \geq T-\Delta$, then $\overline{m}_{T'} \geq \overline{m}_T - \Delta$ holds.
\end{enumerate}
\end{lem}

\begin{proof}
To prove the first assertion set $L:=T \,\mathbf{1}_{D^c}$ and $U:=T \,\mathbf{1}_{D^c} + \mathbf{1}_{D} $. Considering that obviously
$$
\mu_A(\Gamma^{\leq}(L)) \leq \min\big\{\mu_A(\Gamma^{\leq}(T)), \mu_A(\Gamma^{\leq}(T'))\big\} 
\leq \max\big\{\mu_A(\Gamma^{\leq}(T)), \mu_A(\Gamma^{\leq}(T'))\big\} \leq \mu_A(\Gamma^{\leq}(U)) 
$$
as well as $0 \leq \mu_A(\Gamma^{\leq}(U)) - \mu_A(\Gamma^{\leq}(L)) = \mu_A(D \times [0,1])=\lambda(D)$ holds for every $A \in \kc$, the desired inequality follows immediately.\\
To prove the second assertion let $R_\Delta:[0,1] \rightarrow [0,1]$ be defined by $R_\Delta(x)=x+\Delta(mod\, 1)$ and fix $A \in \kc$.
Since obviously $R_\Delta \in \mathcal{T}$, defining $\mu(E \times F)=\mu_A(E \times R_\Delta(F))$ yields a doubly stochastic measure $\mu$ which corresponds to 
a copula $A_\Delta$ (which, in turn, is easily seen to be the transpose of the $R_\Delta$-shuffle $\mathcal{S}_{R_\Delta}(A)$ of $A$). 
Defining $\tilde{T}:[0,1] \rightarrow [0,1]$ by $\tilde{T}(x)=\max\{T(x)-\Delta,0\}$, $\tilde{T} \leq T'$ follows 
and, using disintegration, we get 
\begin{eqnarray*}
\mu_{A_\Delta}(\Gamma^{\leq}(T')) &\geq& \mu_{A_\Delta}(\Gamma^{\leq}(\tilde{T})) = \int_{T^{-1}([\Delta,1])} K_{A_\Delta}\big(x,[0,T(x)-\Delta]\big) d\lambda(x) \\
    &=& \int_{T^{-1}([\Delta,1])} K_{A}\big(x,[\Delta,T(x)]\big) d\lambda(x) \\
    &=& \int_{[0,1]} K_{A}\big(x,[0,T(x)]\big) d\lambda(x) - \int_{T^{-1}([0,\Delta))} K_{A}\big(x,[0,T(x)]\big) d\lambda(x) \\
    && \quad - \,\int_{T^{-1}([\Delta,1])} K_{A}\big(x,[0,\Delta)\big) d\lambda(x) \\ 
    &\geq & \mu_A(\Gamma^{\leq}(T)) - \int_{[0,1]} K_{A}\big(x,[0,\Delta)\big) d\lambda(x) = \mu_A(\Gamma^{\leq}(T)) - \Delta.
\end{eqnarray*}
Since $A \in \kc$ was arbitrary it follows immediately that $\overline{m}_{T'} \geq \overline{m}_T - \Delta$. 
\end{proof}

\noindent 
Slightly modifying the ideas in the first Section of \cite{Ry} it can be shown that for each measurable $T:[0,1] \rightarrow [0,1]$ there exists a non-decreasing function $T^\star:[0,1] \rightarrow [0,1]$ (called the non-decreasing rearrangement of $T$) and a $\lambda$-preserving transformation $\varphi: [0,1] \rightarrow [0,1]$ such that 
\begin{equation} \label{arrangementT}
	T^\star \circ \varphi =T
\end{equation}
holds.
Based on Lemma \ref{lemupper} we can now prove the following main result of this section: 
\begin{thm}\label{main}
Suppose that $T:[0,1]\rightarrow [0,1]$ is measurable. Then we have 
\begin{equation} 
  \overline{m}_T
	= 1+\min_{x \in [0,1]}(x-F_T(x))
  = \overline{m}_{T^\star}.
\end{equation}
\end{thm} 
\begin{proof}
Letting $\mathcal{U}_\varphi:\kc \rightarrow \kc$ denote the operator studied in \cite{TFS} and implicitly defined via
$$
	K_{\mathcal{U}_\varphi(A)}(x,E)=K_{A}(\varphi(x),E),
$$ 
and using disintegration as well as change of coordinates we get that
\begin{eqnarray}\label{Koopman}
	\mu_{\mathcal{U}_\varphi(A)}(\Gamma^{\leq}(T)) &=& \int_{[0,1]} K_{\mathcal{U}_\varphi(A)}(x,[0,T(x)]) d\lambda(x) = 
	\int_{[0,1]} K_{A}\big(\varphi(x),[0,T^\star \circ \varphi(x)]\big) d\lambda(x) \nonumber \\
	&=& \int_{[0,1]} K_A\big(z,[0,T^\star(z)]\big) d\lambda(z) = \mu_A(\Gamma^{\leq}(T^\star))
\end{eqnarray}  
holds for every $A \in \kc$, implying $\overline{m}_T \geq \overline{m}_{T^\star}$. Again using $T^\star \circ \varphi =T$ and the fact that $\varphi$ is $\lambda$-preserving, it is straightforward to verify that $T$ and $T^\star$ have the same d.f., i.e. $F_{T^\star}=F_T$ holds. Therefore, applying Lemma \ref{lemupper} yields
\begin{equation}
	1+\min_{x \in [0,1]}(x-F_{T}(x))=1+\min_{x \in [0,1]}(x-F_{T^\star}(x)) = \overline{m}_{T^\star} \leq \overline{m}_T \leq 1+\min_{x \in [0,1]}(x-F_{T}(x)),
\end{equation} 
from which the desired equality $\overline{m}_{T^\star} = \overline{m}_T$ follows immediately. 
\end{proof}

\begin{figure}[h!]
  \begin{center}
  \includegraphics[width = 10cm,angle=270]{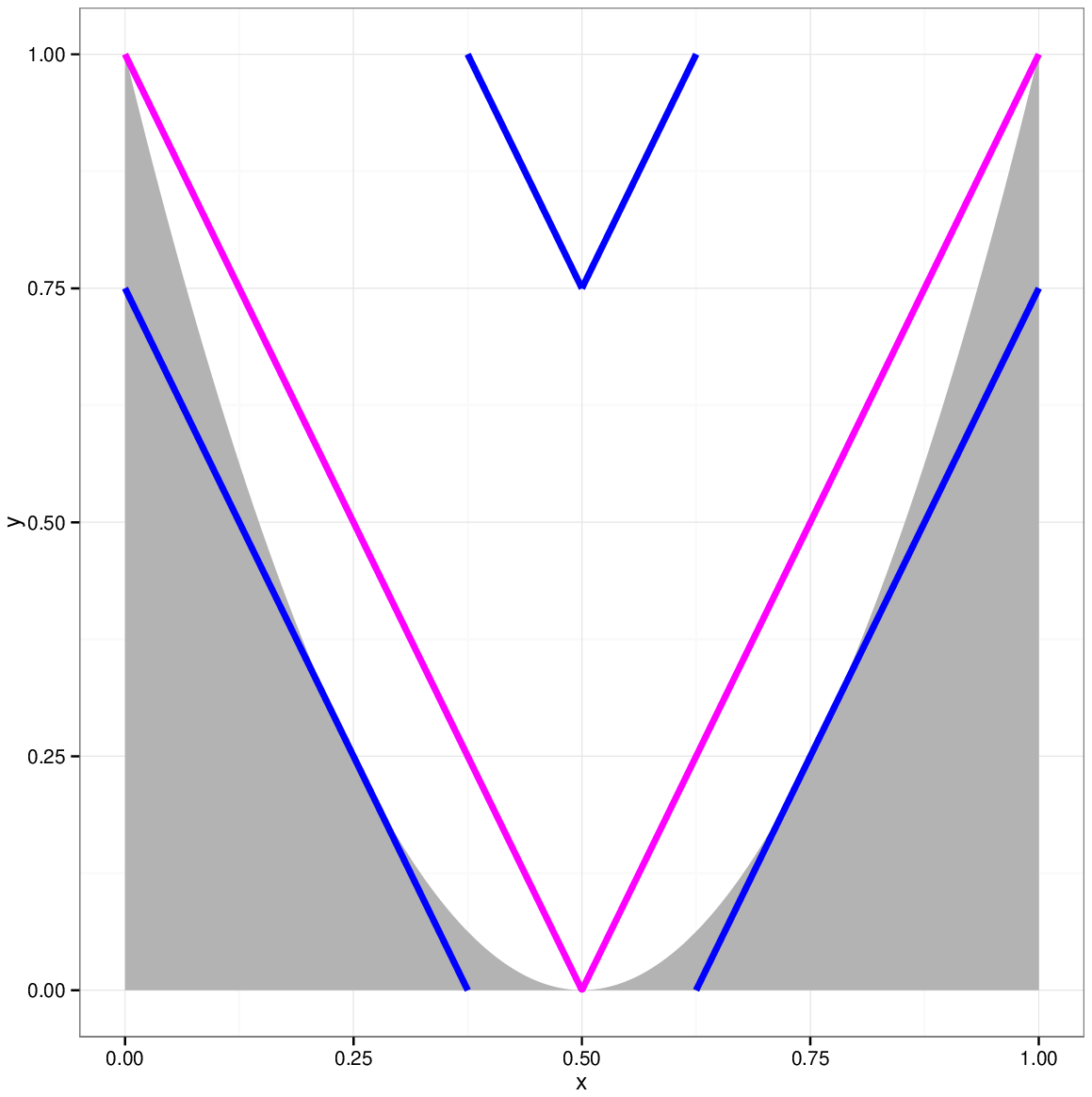}
        \caption{The endograph $\Gamma^\leq(T)$ of the transformation $T$ from Example \ref{examplefinal} (shaded region) as well as the  
        support of the copulas $A_h$ and $A_\varphi$ maximizing/minimizing the mass of $\Gamma^\leq(T)$ (blue and 
        magenta lines, respectively).}\label{exfinal}
  \end{center}
\end{figure}

According to Theorem \ref{solutionmax} the completely dependent copula $A_R \in \kc_d$ fulfills $\overline{m}_{T^\star}=\mu_{A_R}(\Gamma^{\leq}(T^\star))$, so 
eq. (\ref{Koopman}) implies $\mu_{\mathcal{U}_\varphi(A_R)}(\Gamma^{\leq}(T))=\mu_{A_R}(\Gamma^{\leq}(T^\star))=\overline{m}_{T^\star}=\overline{m}_{T}$. 
By defi\-nition of $\mathcal{U}_\varphi(C)$ we have
\begin{equation}\label{cdmax}
K_{\mathcal{U}_\varphi(A_R)}(x,F) = K_{A_R}(\varphi(x),F) = \mathbf{1}_F(R \circ \varphi(x)) =K_{A_{R \circ \varphi}}(x,F),
\end{equation}
so $\mathcal{U}_\varphi(A_R)$ coincides with the completely dependent copula $A_{R \circ \varphi}$
 and the following corollary holds:
\begin{cor}\label{supattained}
Suppose that $T:[0,1]\rightarrow [0,1]$ is measurable. Then there exists a completely dependent copula $A_h \in \kc_d$ such that
$\mu_{A_h}(\Gamma^{\leq}(T))=\overline{m}_T$. 
\end{cor}
Having found a simple analytic formula for the maximal mass of $\Gamma^{\leq}(T)$ we now derive the analogous 
result for the minimal mass and set
\begin{equation}\label{eqdefmin}
\underline{m}_T=\inf_{A \in \kc} \mu_A(\Gamma^{\leq}(T)).
\end{equation}
Given the aforementioned results, the subsequent corollary does not come as a surprise:
\begin{cor}\label{minmass}
For every measurable transformation $T:[0,1]\rightarrow [0,1]$ the following equality holds: 
\begin{equation}\label{eqminmass}
  \underline{m}_T
	= 1-\overline{m}_{1-T}
	= \max_{x \in [0,1]}(x-F_{T}(x-)) 
	= \underline{m}_{T^\star}
\end{equation} 
\end{cor}   
\begin{proof}
We first concentrate on the \emph{strict endograph} $\Gamma^<(T)$, defined by
$$
\Gamma^<(T)=\big\{(x,y) \in [0,1]^2: y<T(x)\big\}.
$$
Defining $T_n:[0,1] \rightarrow [0,1]$ by $T_n(x)=\max\{T(x)-2^{-n},0\}$ for every $x \in [0,1]$ and $n \in \mathbb{N}$ we obviously have that $(\Gamma^{\leq}(T_n))_{n \in \mathbb{N}}$ is 
monotonically increasing and that $\Gamma^<(T)=\bigcup_{n=1}^\infty \Gamma^{\leq}(T_n)$. Lemma \ref{sensitiv} yields $\overline{m}_{T_n} \geq \overline{m}_{T}-2^{-n}$ and 
Corollary \ref{supattained} implies the existence of a copula $A_n \in \kc_d$ with $\mu_{A_n}(\Gamma^{\leq}(T_n))=\overline{m}_{T_n}$. Altogether we get
$$
\overline{m}_{T_n} = \mu_{A_n}(\Gamma^{\leq}(T_n)) \leq \mu_{A_n}(\Gamma^<(T)) \leq \sup_{A \in \kc} \mu_{A}(\Gamma^<(T)) \leq \overline{m}_T,
$$  
so considering $n \rightarrow \infty$ shows that $\sup_{A \in \kc} \mu_{A}(\Gamma^<(T)) = \overline{m}_T$. 
Having this, considering  
\begin{eqnarray*}
\underline{m}_{T} &=& 1 - \sup_{A \in \kc} \mu_A\big( \Gamma^<(1-T)\big) = 1 - \overline{m}_{1-T} = - \min_{x \in [0,1]}(x-F_{1-T}(x)) = \max_{x \in [0,1]}(x-F_{T}(x-)). 
\end{eqnarray*}
yields eq. (\ref{eqminmass}). 
\end{proof}
We close this section with two examples - the first one shows that $\underline{m}_T$ is not necessarily attained whereas the second one focuses on a non-monotonic transformation
for which copulas attaining $\underline{m}_T$ and $\overline{m}_T$ can easily be constructed. 
\begin{ex}
For $T(x)=x$ Corollary \ref{minmass} yields $\underline{m}_{T}=0$. There is, however,
no copula $A$ fulfilling $\mu_A(\Gamma^\leq(T))=0$, i.e. contrary to $\overline{m}_T$, there are situations, in which $\underline{m}_T$ is not attained for any copula. 
Suppose, on the contrary, that $A \in \kc$ fulfills $\mu_A(\Gamma^\leq(T))=0$. Then, defining $h \in \mathcal{T}_b$
by $h(x)=1-x$ and setting $B=\mathcal{U}_h(A)$, we have $\mu_B(\Gamma^\leq(1-T))=0$, so, $B(x,1-x)=0$ holds for every $x\in[0,1]$. The latter implies $B=W$, which is a contradiction since $\mu_W(\Gamma^\leq(1-T))=1$. 
\end{ex}
\begin{ex}\label{examplefinal} 
For $T(x)=4(x-\frac{1}{2})^2$ it is straightforward to find a non-decreasing mapping $T^\star$ and a $\lambda$-preserving
transformation $\varphi$ such that eq. (\ref{arrangementT}) holds. In fact, defining
$\varphi:[0,1] \rightarrow [0,1]$ by
$$
\varphi(x) = \left\{
\begin{array}{rl}
1-2x & \textrm{if } x \in \big[0,\frac{1}{2}\big]\\
-1+2x & \textrm{if } x \in \big(\frac{1}{2},1\big]
\end{array} \right.
$$   
and considering $T^\star(x)=x^2$ we immediately get $T^\star \circ \varphi=T$. Using eq. (\ref{cdmax}), and 
setting $R(x)=x+\frac{3}{4} \,(mod\, 1)$, it follows that
$h=R \circ \varphi $ is $\lambda$-preserving and that $A_h \in \kc_d$ fulfills $\mu_{A_h}(\Gamma^\leq(T))=\overline{m}_T=\overline{m}_{T^\star}=\frac{3}{4}$.  
Considering that for $A_\varphi$ we obviously have $\mu_{A_\varphi}(\Gamma^\leq(T))=1$, we get $\underline{m}_T=0$ which coincides with $\max_{x \in [0,1]}(x-F_{T}(x))$.
Figure \ref{exfinal} depicts the supports of the copulas $A_h$ and $A_\varphi$ as well as the endograph of $T$.  
\end{ex}

\section{Maximizing the mass of the graph and the joint default probability} \label{sec:maxgraph}
In what follows $T:[0,1] \rightarrow [0,1]$ will denote a general non-decreasing transformation.
Since the values at the (at most countably many) discontinuity points of $T$ are irrelevant for the 
maximization problem we will, however, assume that the non-decreasing transformation $T$ 
is right-continuous (or left-continuous if this simplifies technical arguments). 
For every such $T$ there exists a 
set $\Lambda_T \in \mathcal{B}([0,1])$ with $\lambda(\Lambda_T)=1$ such that $T$ is differentiable at every $x \in \Lambda_T$ (see, e.g., \cite{Ru}). In the sequel we will set $T'(x)=0$ for every 
$x \in \Lambda_T^c$ and directly consider $T'$ as integrable function on $[0,1]$ without explicit mentioning. 
Letting $\nu_T$ denote the measure on $\mathcal{B}([0,1])$ generated by 
$T$ via $\nu_T((a,b])=T(b)-T(a)$, it follows that $T'$ is (a version of) the 
Radon-Nikodym derivative of the absolutely continuous component of $\nu_T$ w.r.t. $\lambda$ (see \cite[Chapter 7]{Ru}).
Consequently, for every interval $[a,b] \subseteq [0,1]$ we have 
\begin{equation}\label{inequT'}
\int_{[a,b]} T' d\lambda = \int_{(a,b]} T' d\lambda \leq T(b-)-T(a)=\nu_T((a,b)) \leq \nu_T([a,b]).
\end{equation}  
Inequality (\ref{inequT'}) becomes a chain if equalities for all intervals $(a,b] \subseteq [0,1]$ if and only if $T$
is absolutely continuous.
Define a new measure $\vartheta$ on $\mathcal{B}([0,1])$ by setting 
\begin{equation}\label{eq:vartheta}
\vartheta(E)=\int_{T^{-1}(E)} T' d\lambda.
\end{equation} 
For a given interval $[a,b] \subseteq [0,1]$ we distinguish the following two cases:
(i) If the preimage $T^{-1}([a,b])$ is of the form $[x_1,x_2]$ then using ineq. (\ref{inequT'}) it follows that
\begin{equation*}\label{T'stepI}
\vartheta([a,b])=\int_{T^{-1}([a,b])} T' d\lambda = \int_{[x_1,x_2]} T' d\lambda \leq T(x_2) - T(x_1) \leq b-a=
\lambda([a,b])
\end{equation*}
(ii) If $T^{-1}([a,b])$ is of the form $[x_1,x_2)$ then again by ineq. (\ref{inequT'}) we get
\begin{equation*}
\vartheta([a,b])=\int_{T^{-1}([a,b])} T' d\lambda = \int_{[x_1,x_2)} T' d\lambda \leq T(x_2-) - T(x_1) \leq b-a=
\lambda([a,b]).
\end{equation*}
Having this, the following simple lemma (which will be used in the proof of the main result of this section) is straightforward to prove:
\begin{lem}\label{abscontpi2}
Suppose that $T:[0,1] \rightarrow [0,1]$ is right-continuous and non-decreasing and let $\vartheta$ 
be defined according to eq. (\ref{eq:vartheta}). Then 
$\vartheta(E) \leq \lambda(E)$ holds for every $E \in \mathcal{B}([0,1])$. In particular, $\vartheta$ is absolutely continuous w.r.t. $\lambda$ and the corresponding Radon-Nikodym
derivative $f=\frac{d\vartheta}{d \lambda}$ fulfills $f \leq 1$ $\lambda$-a.e.   
\end{lem}
\begin{proof} Fix $E \in \mathcal{B}([0,1])$ and $\Delta >0 $. By construction of the Lebesgue measure $\lambda$ there exists a family $(I_i)_{i \in \mathbb{N}}$ of compact intervals
fulfilling $E \subseteq \bigcup_{i=1}^\infty I_i$ as well as $\sum_{i=1}^\infty \lambda(I_i) \leq \lambda(E)+\Delta$.
Using $\vartheta([a,b]) \leq \lambda([a,b])$ it follows that 
\begin{eqnarray*}
\vartheta(E)&\leq& \int_{T^{-1}(\bigcup_{i=1}^\infty I_i)} T'd\lambda  \leq \sum_{i=1}^\infty \int_{T^{-1}(I_i)} T'd \lambda \leq \sum_{i=1}^\infty \lambda(I_i) \leq \lambda(E)+\Delta,
\end{eqnarray*}
from which, considering that $\Delta>0$ was arbitrary, we immediately get $\vartheta(E) \leq \lambda(E)$. The remaining assertions are 
straightforward consequences of Radon-Nikodym theorem (\cite{Ru}). 
\end{proof}

As by-product of the results in \cite{DFST} 
we know that for the case of non-singular $T$ (i.e. $\lambda^T$ absolutely continuous w.r.t. $\lambda$) there exists a copula $A \in \kc$ such that, firstly, 
$K_{A}(x,\{Tx\})>0$ for every $x \in [0,1]$ and, secondly, 
$$
\sup_{B \in \kc} \mu_B(\Gamma(T))=\mu_{A}(\Gamma(T))
$$ 
holds. 
If $T$ is not non-singular, there is no copula fulfilling $K_{A}(x,\{Tx\})>0$ for every $x \in [0,1]$ - nevertheless it is possible to find a copula $A^T$ (we write $A^T$ instead of $A_T$ to avoid confusion with completely dependent copulas)
assigning maximal mass
to $\Gamma(T)$ and it is possible to derive a very simple formula for the maximal mass:
\begin{thm}\label{maxmassgraph}
Suppose that $T:[0,1] \rightarrow [0,1]$ is non-decreasing. Then there exists a copula $A^T \in \kc$ such that the following equality holds:
\begin{equation}
\overline{w}_T= \sup_{B \in \kc} \mu_B(\Gamma(T))=\mu_{A^T}(\Gamma(T))=\int_{[0,1]} \min\{T'(x),1\} d\lambda(x)
\end{equation} 
\end{thm}   
\begin{proof} We proceed in several steps and set $a(x)=\min\{T'(x),1\}$ for every $x \in [0,1]$. As first step we show that for every copula $A$ the 
mapping $m_A:[0,1] \rightarrow [0,1]$, defined by $m_A(x)=K_A(x,\{T(x)\})$, fulfills 
$m_A \leq a$ $\lambda$-a.e. Letting $L(m_A)$ denote the set of all Lebesgue points of $m_A$ (see \cite{Ru}) and setting
$\Lambda:=\Lambda_T \cap L(m_A) \cap (0,1)$ it follows that $\lambda(\Lambda)=1$. For every $x \in \Lambda$ and $h>0$ sufficiently small, using disintegration and monotonicity of $T$ we get
\begin{eqnarray}\label{maversusT'}
\frac{1}{2h}\int_{[x-h,x+h]} m_A \,d\lambda &\leq& \frac{1}{2h} \int_{[x-h,x+h]} K_A\big(t,[T(x-h),T(x+h)]\big)\,d\lambda(t)  \\
                                &\leq& \frac{1}{2h}\, \mu_A \big([0,1] \times [T(x-h),T(x+h)] \big) =\frac{T(x+h)-T(x-h)}{2h},\nonumber
\end{eqnarray} 
from which, considering $h \rightarrow0+$ we directly get $0 \leq m_A(x) \leq T'(x)$. Since $x \in \Lambda$ was arbitrary and $m_A(x) \leq 1$ by construction,   
the desired inequality $m_A(x) \leq a(x)$ holds for every $x \in \Lambda$. As direct consequence we get  
$$
\sup_{B \in \kc} \mu_B(\Gamma(T)) \leq \int_{[0,1]} a d\lambda
$$
and the theorem is proved if we can show that there exists a copula $A^T$ fulfilling 
$\mu_{A^T}(\Gamma(T))=\int_{[0,1]} a d\lambda$.  
We distinguish three cases (and, as before, assume w.l.o.g. that $T$ is right-continuous):\\
\textbf{(i)} If $\int_{[0,1]} a\,d\lambda =1$ we get $T' \geq 1$ a.e. 
Since $T'$ is the Radon-Nikodym derivative of the absolutely continuous component of the measure $\nu_T$ mentioned at the beginning of 
this section, considering $\nu_T([0,1])\leq 1$ it follows that $\nu_T$ is absolutely continuous with density $T'$ and that 
$T'=1$ a.e. Hence $\nu_T=\lambda$ and $T=id$, and setting $A^T=M$ yields the desired result $\mu_{A^T}(\Gamma(T))=1$. 
\textbf{(ii)} The case $\int_{[0,1]} a\,d\lambda = 0$ is trivial since every absolutely continuous copula $A$ fulfills $\mu_A(\Gamma(T))=0$.
\textbf{(iii)} In the remaining case of $\int_{[0,1]} a\,d\lambda \in (0,1)$ we can proceed as follows: Define a measure $\mu$ on $\mathcal{B}([0,1]^2)$ by setting
$$
\mu(E \times F) = \int_{E} a(x) \mathbf{1}_F(T(x)) d\lambda(x)
$$    
and extending in the standard way (\cite{Ka,Kl,Ru}) to full $\mathcal{B}([0,1]^2)$. Letting $\pi_1, \pi_2:[0,1]^2 \rightarrow [0,1]$ denote the projections 
onto the first and second coordinate, respectively, for every $E \in \mathcal{B}([0,1])$ we get  
\begin{eqnarray*}
\mu^{\pi_1}(E)&=& \mu(E \times [0,1]) = \int_E a \, d\lambda
\end{eqnarray*}
as well as 
\begin{eqnarray*}
\mu^{\pi_2}(E)&=& \mu([0,1] \times E) = \int_{[0,1]} a(x) \mathbf{1}_E(T(x)) \, d\lambda(x) = \int_{T^{-1}(E)} a \, d\lambda \leq \lambda(E),
\end{eqnarray*}
whereby the last inequality follows from Lemma \ref{abscontpi2}. As direct consequence both $\mu^{\pi_1}$ and $\mu^{\pi_2}$ are absolutely continuous
measures whose densities $f_1,f_2$ fulfill $f_1(x),f_2(x) \in [0,1]$ a.e. and we have $\mu^{\pi_1}([0,1])= \mu^{\pi_2}([0,1]) = \mu([0,1]^2)= \int_{[0,1]} a\, d\lambda \in (0,1)$.
Defining $F_1,F_2:[0,1] \rightarrow [0,1]$ by
$$
F_1(x)=\frac{x-\mu^{\pi_1}([0,x])}{1-\mu^{\pi_1}([0,1])}\, , \quad F_2(x)=\frac{x-\mu^{\pi_2}([0,x])}{1-\mu^{\pi_2}([0,1])}
$$
yields absolutely continuous distribution functions $F_1$ and $F_2$ fulfilling $F_1(0)=F_2(0)=0$. 
Finally, let $R,S:[0,1]^2 \rightarrow [0,1]$ be defined by
\begin{eqnarray}\label{defRS}
R(x_1,x_2)&=& \big(1-\mu^{\pi_1}([0,1])\big) F_1(x_1)F_2(x_2) \\
S(x_1,x_2)&=& \mu\big([0,x_1] \times [0,x_2]\big) \nonumber
\end{eqnarray}
and set $A^T=R+S$. Considering $A^T(1,1)=1$ and the fact that $R$ and $S$ are two-dimensional measure-generating functions by construction, it is now straightforward to
show that $A^T$ is a copula. In fact, the property $A^T(x_1,0)=0$ follows via
$$
A^T(x_1,0)=S(x_1,0) = \mu([0,x_1] \times [0,0]) \leq  \mu([0,1] \times [0,0]) =0 
$$
and the remaining boundary conditions are easily verified too. Since we obviously have $\mu_{A^T}(\Gamma(T))=\mu(\Gamma(T))=\int_{[0,1]} a\, d\lambda$ this completes the proof.
\end{proof}
Notice that in the case of $\int_{[0,1]} \min\{T'(x),1\} d\lambda(x) \in (0,1)$ we could have also defined $R$ by 
$$
R_C(x_1,x_2)= \big(1-\mu^{\pi_1}([0,1])\big) C(F_1(x_1),F_2(x_2)),
$$
whereby $C$ is an arbitrary (not necessarily absolutely continuous) copula, worked with $A_C^T=R_C+S$ and used the fact that in this case 
$\mu_{A_C^T}(\Gamma(T)) \geq \mu(\Gamma(T))=\int_{[0,1]} a\, d\lambda$ holds. As a consequence we get the following corollary:
\begin{cor}\label{formulawt}
If $T:[0,1] \rightarrow [0,1]$ is non-decreasing and $\int_{[0,1]} \min\{T'(x),1\} d\lambda(x) \in (0,1)$, then 
for every copula $C \in \kc$ there exists a copula $A_C^T \in \kc$ such that 
\begin{equation}
\sup_{B \in \kc} \mu_B(\Gamma(T))=\mu_{A_C^T}(\Gamma(T))=\int_{[0,1]} \min\{T'(x),1\} d\lambda(x).
\end{equation} 
\end{cor}
\noindent We now turn to the general problem of calculating 
$$
\overline{w}_T=\sup_{B \in \kc} \mu_B(\Gamma(T))
$$
for general measurable, not necessarily monotonic $T:[0,1] \rightarrow [0,1]$. 
Analogous to the case of $\overline{m}_T$ we first show that rearranging $T$ non decreasingly as $T=T^* \circ \varphi$
does not change the maximum mass, i.e., $\overline{w}_T=\overline{w}_{T^*}$ holds. 
Doing so, we will work with the so-called $\star$-operator (see \cite[Definition 5.4.6]{Du})
$\star :\kc_2 \times \kc_2 \rightarrow \kc_3$, defined by
\begin{equation}
A \star B (x,y,z)=\int_{[0,y]} \partial_{2} A(x,s) \partial_1 B(s,z) d\lambda(s) 
\end{equation}
for all $x,y,z \in [0,1]$. It is straightforward to verify (see \cite{DKQS,Du}) that $\star$ is well-defined, that for all $x,y,z \in [0,1]$ we have 
$A \star B (x,y,1)=A(x,y)$, $A \star B (1,y,z)=B(x,y)$ as well as $A \star B (x,1,z)=A*B(x,z)$, where $*$ denotes the star-product going back to \cite{DNO}.
Furthermore, considering $\partial_{2} A(x,s)=\partial_{1} A^t(s,x)$ it follows immediately that setting
\begin{equation*}\label{kernellift}
K^{13|2}_{A \star B}(y,E \times G) := K_{A^t}(y,E) K_B(y,G)
\end{equation*}
for all $y \in [0,1]$ and $E,G \in \mathcal{B}([0,1])$ and extending in the standard way to $\mathcal{B}([0,1]^2)$ defines a Markov kernel of $A \star B$ w.r.t. the second 
coordinate $y$ (see \cite{FST} and \cite{MFT} for Markov kernels of multivariate copulas). Since $K^{13|2}_{A \star B}(y, \cdot)$ is the product measure of $K_{A^t}(y,\cdot)$ and $K_B(y,\cdot)$, applying 
Fubini's theorem we get that
\begin{equation}\label{fubkern132}
K^{13|2}_{A \star B}(y, \Omega) = \int_{[0,1]} K_B(y,\Omega_x) K_{A^t}(y,dx)
\end{equation}
holds for every $\Omega \in \mathcal{B}([0,1]^2)$.  
\begin{thm}\label{thmrearrange.joint}
Suppose that $T:[0,1] \rightarrow [0,1]$ is measurable and, as before, let $T^*$ with $T^* \circ \varphi =T$ 
denote the non-decreasing rearrangement of $T$. Then $\overline{w}_T=\overline{w}_{T^*}$ holds.
\end{thm}
\begin{proof}
(i) For arbitrary $B \in \kc$, working with $\mathcal{U}_\varphi$, using disintegration and change of coordinates we get
\begin{eqnarray*}
\mu_{\mathcal{U}_\varphi(B)}(\Gamma(T))&=&\int_{[0,1]} K_B\big(\varphi(x),\{T(x)\}\big) d\lambda(x) = \int_{[0,1]} K_B\big(\varphi(x),\{T^*\circ \varphi(x)\}\big) d\lambda(x) \\
&=& \int_{[0,1]} K_B\big(z,\{T^*(z)\}\big) d\lambda(z) = \mu_B(\Gamma(T^*)),
\end{eqnarray*}
from which the inequality $\overline{w}_{T^*} \leq \overline{w}_{T}$ follows immediately. \\
(ii) To prove $\overline{w}_{T^*} \geq \overline{w}_{T}$ we use the $\star$-operator and proceed as follows. 
Letting $A_\varphi $ denote the completely dependent copula induced by $\varphi$, eq. (\ref{fubkern132}) simplifies to
\begin{equation}\label{eqkeylifting}
 K^{13|2}_{A^t_\varphi \star B}(y, \Omega) = \int_{[0,1]} K_B(y,\Omega_x) K_{A_\varphi}(y,dx) =  \int_{[0,1]} K_B(y,\Omega_x) d \delta_{\varphi(y)}(x) = K_B(y,\Omega_{\varphi(y)}).
\end{equation}
Considering $\Omega=\Gamma(T^*)$ we obviously have $\Omega_{\varphi(y)}=\{T^\star\circ\varphi(y)\}=\{T(y)\}$, so it follows that  
$$ 
K^{13|2}_{A^t_\varphi \star B}(y, \Gamma(T^\star)) = K_B(y,\{T(y)\}).
$$
Having this, using disintegration and the fact that $A^t_\varphi \star B (x,1,z)=A^t_\varphi*B(x,z)$ altogether we get
\begin{eqnarray*}
\mu_{A^t_\varphi * B} \big(\Gamma(T^\star) \big) &=& \mu_{A^t_\varphi \star B} \big(\{(x,y,T^\star(x)): x,y \in [0,1]\} \big) = 
       \int_{[0,1]} K^{13|2}_{A^t_\varphi \star B}(y, \Gamma(T^\star)) d\lambda(y) \\
&=& \int_{[0,1]} K_B(y,\{T(y)\}) d\lambda(y) = \mu_B(\Gamma(T)).
\end{eqnarray*}
Considering the fact that $B \in \kc$ was arbitrary the desired inequality $\overline{w}_{T^\star} \geq \overline{w}_{T}$ follows and the theorem is proved. 
\end{proof} 
\begin{rmk}
In the proof of Theorem \ref{thmrearrange.joint} the only properties needed were that $\varphi$ is $\lambda$-preserving and that we have $T^* \circ \varphi =T$ - the fact that $T^*$ is 
non-decreasing was not used. Consequently, 
for arbitrary measurable $S:[0,1] \rightarrow [0,1]$ and arbitrary $\lambda$-preserving $\varphi:[0,1] \rightarrow [0,1]$, setting $T:=S \circ \varphi$ we have
$\overline{w}_{S} = \overline{w}_{T}$.
\end{rmk}
Choosing $C= \mathcal{U}_\varphi(A^{T^*})$, where $A^{T^*}$ denotes the copula maximizing the mass of the graph of the non-decreasing rearrangement
$T^*$ of $T$ directly yields the following result.
\begin{cor}\label{maxmassgraphgeneral}
For every measurable transformation $T:[0,1] \rightarrow [0,1]$ there exists a copula $C \in \kc$ fulfilling 
$\overline{w}_{T} = \mu_C(\Gamma(T))$.
\end{cor}
\noindent Combining Theorem \ref{thmrearrange.joint} and Corollary \ref{formulawt} shows that the identity  
\begin{equation}\label{Tdelta'}
\sup_{B \in \kc} \mu_B(\Gamma(T))=\overline{w}_T= \overline{w}_{T^*}= 
\int_{[0,1]} \min\big\{\big(T^*)'(x),1\big\} d\lambda(x)
\end{equation}
holds for every measurable $T$. In most situations, however, the integral in eq. (\ref{Tdelta'}) is intractable, in particular since calculating the 
rearrangement $T^*$ itself is a nontrivial endeavor. 
Calculating $\overline{m}_T$ for general measurable transformations $T$ in the last section,  
the cumulative distribution function $F_T$ of $T$ plays an important role - we will show now that 
the same is true for $\overline{w}_T$ and derive a very simple formula only involving $F_T$.
\begin{lem}\label{lemcdf}
Suppose that $T:[0,1] \rightarrow [0,1]$ is non-decreasing and let $F_T$ denote the distribution function of $T$. Then $\overline{w}_{T}=\overline{w}_{F_T}$ holds.
\end{lem}
\begin{proof}
 Let $J_T$ denote the set of all discontinuities of $T$ and set 
$I_T=\{y \in [0,1]: \lambda^T(\{y\})>0\}$. Then $I_T$ and $J_T$ are 
at most countably infinite and for every copula $A$ and $y \in I_T$ we have $\mu_A(T^{-1}(\{y\}) \times \{y\})) \leq \mu_A([0,1] \times \{y\})) =0$. 
Setting $N_T:=J_T \cap T^{-1}(I_T)$ obviously $T$ is injective on $N_T$ and for every $A \in \kc$ we have $\mu_A\big((N_T \times [0,1]) \cap \Gamma(T)\big)=0$, implying
$$
\mu_A(\Gamma(T))=\mu_A\big(\underbrace{(N_T^c \times [0,1]) \cap \Gamma(T)}_{=:\Omega_T}\big).
$$
Letting $I_{F_T},J_{F_T},N_{F_T}$ and $\Omega_{F_T}$ denote the corresponding sets for $F_T$ it is straightforward to verify that for every $(x,y) \in [0,1]^2$ we have
$(x,y) \in \Omega_T$ if, and only if $(y,x) \in \Omega_{F_T}$. Having this the desired result follows easily: In fact, letting $A$ denote a copula
with $\mu_{A}(\Gamma(T))=\overline{w}_T$ and considering the transpose $A^t$ we immediately get 
$$
\overline{w}_T=\mu_{A}(\Gamma(T))=\mu_{A}(\Omega_T)=\mu_{A^t}(Q_{F_T}) = \mu_{A^t}(\Gamma(F_T)) 
\leq \overline{w}_{F_T}. 
$$
Since the other inequality follows in the same manner the desired equality is proved. 
\end{proof}  
\noindent Considering that $T$ and $T^*$ have the same distribution function and applying Lemma \ref{lemcdf} yields a handier version of eq. (\ref{Tdelta'}):
\begin{cor}
For every measurable $T:[0,1] \rightarrow [0,1]$ the following equality holds: 
\begin{equation}\label{formulawtfinal}
\sup_{B \in \kc} \mu_B(\Gamma(T))=\overline{w}_T= \overline{w}_{F_T}=\int_{[0,1]} \min\big\{F_T'(x),1\big\} d\lambda(x)
\end{equation}
\end{cor} 

\noindent Notice that in case $T:[0,1] \rightarrow [0,1]$ is non-decreasing and continuous and fulfills 
$T'=0$ $\lambda$-almost everywhere according to eq. (\ref{formulawt}) $\overline{w}_T=0$, i.e., no copula
assigns mass to $\Gamma(T)$. This result  is not surprising - considering the fact, however, that in the 
language of Baire categories a `typical' monotonic function is singular (as established in \cite{Zam}) 
we could infer that copulas assign no mass to `typical' monotonic functions, which seems quite counterintuitive.  
In \cite[Theorem 3]{DFST} it was shown that for every non-singular $T:[0,1] \rightarrow [0,1]$ there exists a copula $A$ such that the singular component of $A$
is concentrated on $\Gamma(T)$ and that we have $K_A(x, \{T(x)\}) > 0$ for $\lambda$-almost every $x \in [0,1]$. 
Based on Theorem \ref{maxmassgraph} and Theorem \ref{thmrearrange.joint} we can give a necessary and sufficient condition for the existence of a copula $A$ 
fulfilling $K_A(x, \{T(x)\}) > 0$ for every $x \in [0,1]$ in terms of the non-decreasing rearrangement $T^*$ of 
$T$ and in terms of absolute continuity of $F_T$. 
\begin{cor}\label{masseverywhere}
Suppose that $T:[0,1] \rightarrow [0,1]$ is measurable, let $T^*$ denote its non-decreasing rearrangement and $F_T$ its distribution function. Then the following conditions are
equivalent:
\begin{enumerate}
\item[(a)] There exists a copula $C$ fulfilling $K_C(x, \{T(x)\}) > 0$ for $\lambda$-almost every $x \in [0,1]$.
\item[(b)] $(T^*)'(x) >0$ for $\lambda$-almost every $x \in [0,1]$.
\item[(c)] $\lambda^T$ is absolutely continuous.
\item[(d)] $F_T$ is absolutely continuous.
\end{enumerate}
\end{cor}
\begin{proof}
It is clear that (c) and (d) are equivalent so it suffices to prove $(b) \Rightarrow (a) \Rightarrow (c) \Rightarrow (b)$, which can be done as follows.
(i) Let $A$ be a copula assigning maximum mass to $\Gamma(T^*)$ as constructed in the proof of Theorem \ref{maxmassgraph}. 
Letting $K_{A}(\cdot,\cdot)$ denote a version of the Markov kernel of $A$ fulfilling $K_{A}(x,\{T^*(x)\})>0$ 
for $\lambda$-almost every $x \in [0,1]$ and considering $C=\mathcal{U}_\varphi(A)$ we get
\begin{equation} \label{formula.characterization}
K_C(x,\{T(x)\})= K_{A}(\varphi(x),\{T(x)\}) = K_{A}(\varphi(x),\{T^* \circ \varphi(x)\}) >0,
\end{equation}
so (b) implies (a). \\
The implication $(a) \Rightarrow (c)$ is a direct consequence of the fact that for every $N \in \mathcal{B}([0,1])$ with $\lambda(N)=0$ 
and an arbitrary copula $C$ fulfilling (a) we have
\begin{eqnarray*}
0&=& \lambda(N)= \mu_C\big([0,1] \times N \big) \geq  \mu_C\big(T^{-1}(N) \times N \big) \geq \int_{T^{-1}(N)} \underbrace{K_C(x,\{T(x)\})}_{>0} d\lambda(x),
\end{eqnarray*}
from which $\lambda^T(N)=0$ follows immediately. \\
(iii) Simplifying notation set $S:=T^*$ and suppose now that $\lambda^{S}=\lambda^T$ is absolutely continuous. We want to show that $S'(x)>0$ for 
$\lambda$-almost every $x \in [0,1]$. Since $S$ is non-decreasing, 
considering that $\lambda^S(\{y\})=0$ and that $S^{-1}(\{y\})$ is an interval for every $y \in [0,1]$, it follows that $S^{-1}(\{y\})$ is either empty or a degenerated
interval consisting of one single point, so $S$ is necessarily strictly increasing on $[0,1]$. Additionally, for every $E \subseteq [0,1]$ we obviously have 
$S^{-1} (S (E))=E$. Assume that $S(0)=0$ (if $S(0)>0$ holds proceed with the function $\tilde{S}$ that coincides with $S$ on $(0,1]$ and fulfills $\tilde{S}(0)=0$).
Letting $f$ denote the Radon-Nikodym derivative of $\lambda^{S}$ w.r.t. $\lambda$
we may w.l.o.g. assume $0 \leq f(z)<\infty$ for every $z \in [0,1]$. The function $g:[0,1] \rightarrow [0,1]$, defined by
$y \mapsto \int_{[0,y]} f d\lambda$ is non-decreasing and
fulfills 
\begin{equation}\label{eqgTDelta}
g(S(x))=\lambda^{S}([S(0),S(x)])=\lambda^{S}(S([0,x])=x.
\end{equation}
Choose $\Psi,\Lambda \in \mathcal{B}([0,1])$ 
with $\lambda(\Lambda)=1=\lambda(\Psi)$ in such a way that $g$ is differentiable at every $y \in \Lambda$ and fulfills $g'(y)=f(y)$ and that  $S$ is differentiable at every $z \in \Psi$. 
For every $x \in S^{-1}(\Lambda) \cap \Psi$ applying the chain rule together with equ. (\ref{eqgTDelta})
yields 
$$
1=f(S(x)) S'(x),
$$ 
hence $S'(x)>0$. This completes the proof since $\lambda(S^{-1}(\Lambda) \cap \Psi)=1$.  
\end{proof}

\begin{rmk}
Again using the $\star$-operator allows for a direct proof of the implication $(a) \Rightarrow (b)$ of Corollary \ref{masseverywhere}: Suppose that $E \in \mathcal{B}([0,1])$ is arbitrary but fixed. Applying eq. (\ref{eqkeylifting}) to the set 
$\Gamma_E(T^*):=\{(x,T^*(x)): x \in E\}=\Gamma(T^*) \cap (E \times [0,1])$ for every $B \in \kc$ we get 
$$
K^{13|2}_{A^t_\varphi \star B}(y, \Gamma_E(T^*)) = \mathbf{1}_{E}(\varphi(y)) \,K_B(y,\{T(y)\}),
$$
so, using disintegration 
\begin{eqnarray*}
\mu_{A^t_\varphi * B} \big(\Gamma_E(T^*) \big) &=& \mu_{A^t_\varphi \star B} \big(\{(x,y,T^*(x)): x \in E,y \in [0,1]\} \big) = 
       \int_{[0,1]} K^{13|2}_{A^t_\varphi \star B}(y, \Gamma_E(T^*)) d\lambda(y) \\
&=& \int_{[0,1]} \mathbf{1}_{E}(\varphi(y)) \,K_B(y,\{T(y)\}) d\lambda(y) = \int_{\varphi^{-1}(E)} \,K_B(y,\{T(y)\}) d\lambda(y)
\end{eqnarray*}
follows. Suppose now that $B \in \kc$ fulfills $K_B(y,\{T(y)\})>0$ for every $y \in [0,1]$. Using the fact that $\varphi$ is $\lambda$-preserving we get that
$\mu_{A^t_\varphi * B} \big(\Gamma_E(T^*) \big)> 0$ if, and only if $\lambda(E)>0$. 
Since for $E:=\{x \in [0,1]: K_{A^t_\varphi * B}(x,\{T^*(x)\})=0\} \in \mathcal{B}([0,1])$
obviously $\mu_{A^t_\varphi * B} \big(\Gamma_E(T^*) \big)=0$ holds, the latter implies $\lambda(E)=0$, so we can find a version of the kernel $K_C(\cdot,\cdot)$ of 
the copula $C=A^t_\varphi * B$ such that $K_C(x,\{T^\star(x)\})>0$ holds for $\lambda$-almost every $x \in [0,1]$. Having this, proceeding analogously to (\ref{formula.characterization}) directly yields condition (a). $\blacksquare$
\end{rmk}


\section{When $\overline{m}_T$ and $\overline{w}_T$ coincide}\label{sec:coincide}
The results in the previous two sections allows to characterize all 
non-decreasing functions $T:[0,1] \rightarrow [0,1]$ for which the maximum mass of 
the graph $\Gamma(T)$ and the maximum mass of the endograph $\Gamma^\leq(T)$ coincide: 
\begin{thm}\label{maxmassconcide}
Suppose that $T:[0,1] \rightarrow [0,1]$ is non-decreasing and let $\Lambda_T$ denote the set of all points at which $T$ is differentiable. 
Then the following two assertions are equivalent.
\begin{enumerate}
\item[(a)] $\overline{m}_T = \overline{w}_T$.
\item[(b)] $T(0)=0$ and there exists a point $x_0 \in [0,1]$ such that the following conditions hold: 
\begin{itemize}
\item[(i)] $T$ is absolutely continuous on $[0,x_0]$, 
\item[(ii)] $\Omega_0:=\{x \in [0,x_0] \cap \Lambda_T:\, T'(x) \leq 1 \}$ fulfills $\lambda(\Omega_0)=x_0$, 
\item[(iii)] $\Omega_1:=\{x \in [x_0,1] \cap \Lambda_T:\, T'(x) \geq 1 \}$ fulfills $\lambda(\Omega_1)=1-x_0$.
\end{itemize}
\end{enumerate}
\end{thm}
\begin{proof}
We may, w.l.o.g., assume that $T$ is left-continuous. \\
\textbf{(I)} Suppose that $T$ fulfills the second assertion. It follows immediately from 
condition (i) that for every $z \in [0,x_0]$ we have $T(z)=\int_{[0,z]} T' d\lambda$, hence, by condition (ii), the mapping $z \mapsto T(z)-z=\int_{[0,z]} (T'-1) d\lambda$ 
is non-increasing on $[0,x_0]$ and we have $\inf_{z \in [0,x_0]} (T(z)-z) = T(x_0)-x_0$. Additionally, condition (iii) implies that $z \mapsto T(z)-z$ is non-decreasing on
$[x_0,1]$, from which, using Theorem \ref{solutionmax} we altogether get  
$$
\sup_{B \in \kc} \mu_B(\Gamma^\leq(T)) = \overline{m}_T= 1+\inf_{x \in [0,1]} (T(x)-x)=1+T(x_0)-x_0.
$$
Taking into account that (i)-(iii) also imply 
$$
\int_{[0,1]}  \min\{T'(x),1\} d\lambda(x) =  \int_{[0,x_0]} T' d\lambda + \int_{[x_0,1]} 1 d\lambda = T(x_0)+1-x_0
$$
the desired equality $\overline{m}_T = \overline{w}_T$ follows. \\
\textbf{(II)} On the other hand, if $\overline{m}_T = \overline{w}_T$ holds, then again by Theorem \ref{solutionmax}, 
left-continuity of $T$ and Theorem \ref{maxmassgraph}, there exists some $x_0 \in [0,1]$ such that 
$$
1+T(x_0)-x_0 =  \overline{m}_T = \overline{w}_T = \int_{[0,1]} \min\{T'(x),1\} d\lambda(x)
$$
holds. Considering $\int_{[0,x_0]} \min\{T'(x),1\} d\lambda(x) \leq T(x_0)-T(0) \leq T(x_0)$ together with the fact that 
$\int_{[x_0,1]} \min\{T'(x),1\} d\lambda(x) \leq 1-x_0$ it follows immediately that $T$ has to fulfill $T(0)=0$ as 
well as
\begin{eqnarray*}
\int_{[0,x_0]} \min\{T'(x),1\} d\lambda(x)&=&T(x_0), \quad \int_{[x_0,1]} \min\{T'(x),1\} d\lambda(x)= 1-x_0. 
\end{eqnarray*} 
The latter, however, implies that $T$ is absolutely continuous on $[0,x_0]$ and that $T$ fulfills (ii) and (iii). 
\end{proof}
\noindent Considering that every convex function $T:[0,1] \rightarrow [0,1]$ with $T(0)=0$ fulfills the properties listed in condition (b) of Theorem \ref{maxmassconcide} we immediately get the following result:
\begin{cor}\label{convexT}
If $T:[0,1] \rightarrow [0,1]$ is convex and fulfills $T(0)=0$ then $\overline{m}_T = \overline{w}_T$ holds.
\end{cor}
\noindent According to Corollary \ref{convexT}, given a non-decreasing transformation $T:[0,1] \rightarrow [0,1]$, convexity and $T(0)=0$ is sufficient 
for $\overline{m}_T=\overline{w}_T$. The two conditions are, however, far from being necessary - the following 
example shows that equality can also hold for non-decreasing transformations that are not even locally convex.
\begin{ex}
Let $\Omega \in \mathcal{B}([0,1])$ denote a set with $\lambda(\Omega)=\frac{1}{2}$ such that $\lambda((a,b) \cap \Omega)>0$ and $\lambda((a,b) \cap \Omega^c)>0$ hold 
for every non-empty open interval $(a,b) \subseteq [0,1]$ (for a possible construction see \cite[Lemma 3.1]{FSTComm}). 
Define the function $S: [0,1] \rightarrow [0,1]$ by $S(x)=\int_{[0,x]} \mathbf{1}_\Omega(y) \,d \lambda(y)$ for every $x \in [0,1]$. 
Then $S$ is strictly increasing, $S(0)=0$, $S(1)=\frac{1}{2}$, $S$ is absolutely continuous and 
$S'(x)=\mathbf{1}_\Omega(x) \in \{0,1\}$ holds for $\lambda$-almost every $x \in [0,1]$ (see \cite{Ru}). There exists a unique $x_0 \in (\frac{1}{2},1)$ fulfilling $S(x_0)=1-x_0 < \frac{1}{2}$ and 
the properties of $\Omega$ imply that $S$ is not convex on any non-degenerated subinterval of $[0,1]$. Based on $S$ define a new transformation $T:[0,1] \rightarrow [0,1]$ by
$$
T(x) = \left\{
\begin{array}{rl}
S(x) & \textrm{if } x \in [0, x_0],\\
\frac{x-x_0}{2(1-x_0)} + S(x) & \textrm{if } x \in [x_0,1].
\end{array} \right.
$$
It is straightforward to verify that $T:[0,1] \rightarrow [0,1]$ is a strictly increasing homeomorphism of $[0,1]$,
which fulfills all properties stated in assertion (b) of Theorem \ref{maxmassconcide}. $T$ is, however,  
obviously not convex on any non-degenerated subinterval of $[0,1]$ (since its derivative is not non-decreasing on any
non-empty open interval). 
\end{ex} 
\begin{rmk}
Let $T:[0,1] \rightarrow [0,1] $ be non-decreasing and right-continuous, and assume that $T(0)=0$ holds. 
According to Theorem \ref{maxmassconcide} in order to have $ \overline{m}_T = \overline{w}_T$ 
the transformation $T$ needs to be absolutely continuous on the interval $[0,x_0]$ - on the interval $[x_0,1]$, 
however, $T$ (interpreted as univariate measure-generating function) may be also have a non-degenerated
discrete and/or singular component on $[x_0,1]$ as long as $T' \geq 1$ holds $\lambda$-almost everywhere on $[x_0,1]$.  
\end{rmk}
\begin{rmk}
The proof of Theorem \ref{maxmassconcide} also shows that for non-decreasing $T:[0,1] \rightarrow [0,1]$ 
all copulas $A$ with $\mu_A(\Gamma(T))=\overline{m}_T = \overline{w}_T$ fulfill 
the following three conditions:
\begin{eqnarray*}
\mu_A\big([0,x_0]\times [0,T(x_0)]\big) &=& \mu_A\big(([0,x_0]\times [0,1]) \cap \Gamma(T)\big)=T(x_0) \\
\mu_A\big([x_0,1]\times [T(x_0),1]\big) &=&  \mu_A\big(([x_0,1]\times [T(x_0),1]) \cap \Gamma(T)\big)=1-x_0 \\
\mu_A\big(([x_0,1]\times [0,T(x_0)]) &=& 0
\end{eqnarray*} 
\end{rmk}
\noindent Again working with non-decreasing rearrangements and using the previous results yields 
the following corollary:
\begin{cor}\label{graphendoviarearrangement}
For every measurable $T:[0,1] \rightarrow [0,1]$ the following two conditions are equivalent (as before $T^*$ denotes the
non-decreasing rearrangement and $\Lambda_{T^*} \in \mathcal{B}([0,1])$ the set of all points at which
$T^*$ is differentiable):
\begin{itemize}
\item[(a)] $\overline{m}_T = \overline{w}_T$.
\item[(b)]  $T^*(0)=0$ and there exists some 
$x_0 \in [0,1]$ such that the following two properties hold: 
\begin{itemize}
\item [(i)] $T^*$ is absolutely continuous on $[0,x_0]$, 
\item [(ii)] $\Omega_0:=\{x \in [0,x_0] \cap \Lambda_{T^*}: (T^\star)'(x) \leq 1 \}$ fulfills $\lambda(\Omega_0)=x_0$, 
\item [(iii)] $\Omega_1:=\{x \in [x_0,1] \cap \Lambda_{T^*}: (T^\star)'(x) \geq 1 \}$ fulfills $\lambda(\Omega_1)=1-x_0$.
\end{itemize}
\end{itemize}
\end{cor}


\section{Estimating the maximum probability of a prior default} \label{sec:asymptotics}
Throughout this section we assume that $F$ and $G$ are univariate continuous distribution functions, let $T$ be defined 
by $T:=G \circ F^-$ on $(0,1)$ and set $T(0)=0$ and $T(1)=T(1-)$, which implies that $T$ is left-continuous on $[0,1]$. 
Notice that for such $T$ there exists some (not necessarily unique) $x \in (0,1]$ with $\overline{m}_T=1+T(x)-x$.  
If $X_1,\ldots,X_n$ and $Y_1,\ldots,Y_n$ are independent samples of $F$ and $G$, respectively, then 
it seems natural to estimate $\overline{m}_T$ by $\overline{m}_{T_n}$ where $T_n=G_n \circ F_n^-$ and $F_n,G_n$ are the 
empirical distribution functions corresponding to $X_1,\ldots,X_n$ and $Y_1,\ldots,Y_n$. 
We are now going to show that $\overline{m}_{T_n}$ is a strongly consistent estimator for $\overline{m}_{T}$ and 
start with the following simple lemma.
\begin{lem}\label{lem:est}
Suppose that $F$ and $G$ are continuous univariate distribution functions. Then with probability one $\lim_{n \rightarrow \infty} \vert T_n(u)-T(u) \vert=0$ holds for every continuity point 
$u \in (0,1)$ of $F^-$. In particular $(T_n)_{n \rightarrow \infty}$ converges to $T$ $\lambda$-almost everywhere.  
\end{lem}
\begin{proof}
Glivenko-Cantelli theorem implies that with probability we have uniform convergence of $(F_n)_{n \in \mathbb{N}}$ to $F$ 
and of $(G_n)_{n \in \mathbb{N}}$ to $G$. Applying Lemma 21.2 in \cite{vdV} it follows that for every
continuity point $u \in (0,1)$ of $F^-$ we have $\lim_{n \rightarrow \infty} F_n^-(u)=F^-(u)$ from which the desired
result follows by a straightforward application of the triangle inequality.
\end{proof}

\begin{thm}
Suppose that $F$ and $G$ are continuous distribution functions and let $X_1,\ldots,X_n$ and $Y_1,\ldots,Y_n$ be 
independent samples of $F$ and $G$, respectively. Then with probability one we have 
$\lim_{n \rightarrow \infty} \overline{m}_{T_n}=\overline{m}_T$, i.e. $\overline{m}_{T_n}$ is a strongly consistent estimator of $\overline{m}_T$.
\end{thm}
\begin{proof}
According to Lemma \ref{lem:est} we may assume that $(T_n)_{n \in \mathbb{N}}$ converges to $T$ $\lambda$-almost everywhere.
(i) Fix $\varepsilon >0$ and suppose that $x \in (0,1]$ fulfills $\overline{m}_T=1+T(x)-x$. Then there exists some
$z \in (x-\varepsilon,x) $ such that $z$ is a continuity point of $F^-$ and according to Lemma \ref{lem:est} we can 
find an index $n_0 \in \mathbb{N}$ such that $\vert T_n(z)-T(z)\vert < \varepsilon$, hence 
\begin{eqnarray*}
\overline{m}_{T_n} &\leq&  1 + T_n(z) - z \leq  1 + T(z) + \varepsilon - z \leq 1 + T(x) + \varepsilon - z < 1 + T(x) + \varepsilon - x + \varepsilon  \\
&=& 1+ T(x)-x + 2 \varepsilon = \overline{m}_T - 2 \varepsilon
\end{eqnarray*}
for every $n \geq n_0$. Considering that $\varepsilon>0$ was arbitrary 
$\limsup_{n \rightarrow \infty} \overline{m}_{T_n} \leq \overline{m}_T$ follows.\\
(ii) Suppose now that $\liminf_{n \rightarrow \infty} \overline{m}_{T_n} = \overline{m}_T - 2 \Delta$ holds for some 
$\Delta>0$. Without loss of generality (choose an appropriate subsequence if necessary) we may assume that  
$$
\lim_{n \rightarrow \infty} \overline{m}_{T_{n}} = \overline{m}_T - 2 \Delta.
$$   
Then for every $n \in \mathbb{N}$ there exists some $x_n \in (0,1]$ with $1+T_n(x_n) - x_n < \overline{m}_{T_{n}}
 + \frac{\Delta}{2}$ and we can find an index $n_0$ such that for every $n \geq n_0$ we have 
$ \overline{m}_{T_{n}} <  \overline{m}_{T} - \frac{3 \Delta}{2}$ and  
\begin{eqnarray*}
1+T_n(x_n) - x_n < \overline{m}_{T_{n}}  + \tfrac{\Delta}{2} <  \overline{m}_{T} - \Delta.
\end{eqnarray*} 
Compactness of $[0,1]$ implies the existence of a subsequence $(x_{n_j})_{j \in \mathbb{N}}$ with limit $x \in [0,1]$.
Now, choose $\delta \in (0,\frac{\Delta}{4})$ so that $x-\delta$ is a continuity point of $F^-$ and
 $T(x-\delta) \geq T(x)-\frac{\Delta}{4}$ holds. 
Choose $j_0 \in \mathbb{N}$ in such a way that $n_{j_0}\geq n_0$ and that 
$\vert x_{n_j} - x \vert  \leq \delta$ for every $j \geq j_0$. 
According to Lemma \ref{lem:est} we can find another index $j_1 \in \mathbb{N}$ in such a way that $n_{j_1}\geq n_{j_0}$ and that
$\vert T_{n_j}(x-\delta) - T(x-\delta) \vert \leq \frac{\Delta}{4}$ for every $j \geq j_1$. 
Then for $j \geq j_1$ we altogether get
\begin{eqnarray*}
  1+T(x)-x 
	&\leq& 1 + T(x-\delta) + \tfrac{\Delta}{4} -x 
	 \leq  1 + T_{n_j}(x-\delta) + \tfrac{\Delta}{2} -x 
	\\
	&\leq& 1 + T_{n_j}(x_{n_j}) + \tfrac{\Delta}{2} -x 
   \leq  1 + T_{n_j}(x_{n_j}) + \tfrac{\Delta}{2} -x_{n_j} + \tfrac{\Delta}{4} 
	\\
	&  < & \overline{m}_T - \Delta +  \tfrac{3\, \Delta}{4} 
= \overline{m}_T -  \tfrac{\Delta}{4},
\end{eqnarray*}
a contradiction to the definition of $\overline{m}_T$. This shows
$\liminf_{n \rightarrow \infty} \overline{m}_{T_n} \geq \overline{m}_T$ and the proof is complete.
\end{proof}
As final step we will now show that under mild regularity conditions on $T$ (or, equivalently on $F$ and $G$) 
the estimator $\overline{m}_{T_n}$ is asymptotically normal. 
To derive asymptotic normality we will apply the functional Delta method (see \cite{vdV}) and build upon the following
two lemmata, whereby as in \cite{vdV} $[a,b] \subseteq [-\infty,\infty]$ and $\mathbb{D}([a,b])$ 
will denote the family of all cadlag functions endowed with the uniform distance $\Vert \cdot \Vert_\infty$: 
\begin{lem}
\label{lem:Gx-had-diff}
Define $\phi\colon [a,b] \times \mathbb{D}[a,b] \rightarrow \mathbb{R}$ by $\phi(x,G) = G(x)$ and suppose that 
$G \in \mathbb{D}[a,b]$ is differentiable at $x \in (a,b)$. Then $\phi$ is Hadamard differentiable at $(x,G)$ tangentially to 
the set of tuples $(h_1,h_2) \in \mathbb{R} \times \mathbb{D}[a,b]$ where $h_2$ is continuous at $x$, with derivative 
$\phi^\prime: \mathbb{R} \times \mathbb{D}[a,b] \to \mathbb{R}$ fulfilling
$\phi^\prime(h_1,h_2) = h_1 G^\prime(x) + h_2(x)$.
\end{lem}
\begin{proof}
Let $h_{1,t} \rightarrow h_1$, $h_{2,t} \rightarrow h_2$ and $t \rightarrow 0$
such that $x+th_{1,t} \in [a,b]$ for sufficiently small $t$. 
Then using Taylor's formula we get
\begin{align*}
	\frac{\phi(x+th_{1,t},G+th_{2,t})-\phi(x,G)}{t} &= \frac{(G+th_{2,t})(x+th_{1,t})-G(x)}{t}\\
	&= \frac{G(x+th_{1,t})+th_{2,t}(x+th_{1,t})-G(x)}{t}\\
	&= \frac{G(x)+G^\prime(x)th_{1,t}+o(t)-G(x)}{t} + h_{2,t}(x+th_{1,t})\\
	&\rightarrow G^\prime(x)h_{1} + h_{2}(x),
\end{align*}
where in the last step we used continuity of $h_2$ at $x$. 
\end{proof} 
\begin{lem} \label{lem:Tx-had-diff}
Define $\phi\colon \mathbb{D}[a,b] \times \mathbb{D}[a,b] \rightarrow \mathbb{R}$ by $\phi(F,G) = 
G\circ F^{-}(p)$, consider $p \in (0,1)$ and set $x_p = F^{-}(p) \in (a,b)$. 
Furthermore let $F, G \in \mathbb{D}[a,b]$ be differentiable at $x_p$ with
$F'(x_p)>0$. 
Then $\phi$ is Hadamard differentiable at $(F,G)$ tangentially to $(h_1,h_2) \in \mathbb{D}[a,b] \times \mathbb{D}[a,b]$ where $h_1$ and $h_2$ are continuous at $x_p$, with derivative $\phi^\prime: \mathbb{D}[a,b] \times \mathbb{D}[a,b] \mapsto \mathbb{R}$ fulfilling
$$
\phi^\prime(h_1,h_2) = - \,h_1(x_p)\,\frac{G^\prime(x_p)}{F^\prime(x_p)} + h_2(x_p).
$$
\end{lem}
\begin{proof} 
As a consequence of \cite[Lemma 21.3]{vdV}, 
the map $\phi_1\colon \mathbb{D}[a,b] \times \mathbb{D}[a,b] \rightarrow \mathbb{R} \times \mathbb{D}[a,b]$ defined by $\phi_1(F,G) = (F^{-}(p),G)$ 
is Hadamard differentiable at $(F,G)$ tangentially to the set of functions 
$(h_1,h_2) \in \mathbb{D}[a,b] \times \mathbb{D}[a,b]$ where $h_1$ is continuous at $x_p$, 
with derivative $\phi_1^\prime: \mathbb{D}[a,b] \times \mathbb{D}[a,b] \to \mathbb{R} \times \mathbb{D}[a,b]$ fulfilling
$\phi_1^\prime(h_1,h_2) = (-h_1(x_p)/F^\prime(x_p),h_2)$.
According to Lemma \ref{lem:Gx-had-diff} the map $\phi_2\colon [a,b] \times \mathbb{D}[a,b] \rightarrow \mathbb{R}$ 
defined by $\phi_2(x,G) = G(x)$ is Hadamard differentiable at $(F^{-}(p),G)$ tangentially to 
the set of tuples $(h_1,h_2) \in \mathbb{R} \times \mathbb{D}[a,b]$ where $h_2$ is continuous at $x_p$, with derivative 
$\phi^\prime: \mathbb{R} \times \mathbb{D}[a,b] \to \mathbb{R}$ fulfilling
$\phi^\prime(h_1,h_2) = h_1 G^\prime(x) + h_2(x)$.
It hence follows from the Chain rule for Hadamard derivatives (see \cite[Theorem 20.9]{vdV}) that the 
transformation $\phi_2\circ\phi_1$ is Hadamard differentiable as well, which completes the proof. 
\end{proof}

Given an interval $[a,b] \subseteq \mathbb{R}$, let $\mathbb{D}_1$ denote the set of all restrictions of distribution functions on $\mathbb{R}$ to $[a,b]$,
and $\mathbb{D}_2$ the subset of $\mathbb{D}_1$ consisting of all distribution functions of probability measures assigning mass $1$ to $(a,b]$.
Furthermore let $\mathbb{C}[a,b]$ denote the family of all continuous functions on $[a,b]$.
The following corollary works analogously to Lemma 21.4 in \cite{vdV}.
\begin{cor} \label{cor:FG-had-diff} 
\begin{enumerate} 
\item 
Let $0<p_1<p_2<1$ and let $F,G$ be continuously differentiable on the interval $[a,b] = [F^{-}(p_1)-\epsilon,F^{-}(p_2)+\epsilon]$ for some $\epsilon > 0$, with the derivative of $F$ being strictly positive. 
Then $\phi\colon \mathbb{D}_1 \times \mathbb{D}[a,b] \rightarrow \mathbb{D}[0,1]$ defined by $\phi(\overline{F},\overline{G}) = \overline{G}\circ \overline{F}^{-}$ is Hadamard differentiable at $(F,G)$ tangentially to $\mathbb{C}[a,b] \times \mathbb{C}[a,b]$.
\item 
Let $F$ have compact support $[a,b]$ and let $F,G$ be continuously differentiable on $[a,b]$ with the derivative of $F$ being strictly positive. 
Then $\phi\colon \mathbb{D}_2 \times \mathbb{D}[a,b] \rightarrow \mathbb{D}[0,1]$ defined by $\phi(\overline{F},\overline{G}) = \overline{G}\circ \overline{F}^{-}$ is Hadamard differentiable at $(F,G)$ tangentially to $\mathbb{C}[a,b] \times \mathbb{C}[a,b]$.
\end{enumerate}
In both cases the derivative is the map $$(h_1,h_2) \mapsto \left(-\,h_1\,\frac{G^\prime}{F^\prime} + h_2\right) \circ F^{-}.$$
\end{cor}
The next result is immediate from \cite{CCR}:

\begin{lem} \label{lem:m-had-diff}
Define $\phi\colon \mathbb{D}[0,1] \rightarrow \mathbb{R}$ as $\phi(T) = 1 + \inf_{x\in[0,1]} (T(x) - x)$. 
Let $T \in \mathbb{D}[0,1]$ be such that there exists a unique $x^*\in (0,1)$ with $1 + T(x^*-) - x^* = 1 + \inf_{x\in[0,1]} T(x) - x$. 
Then $\phi$ is Hadamard differentiable at $T$ tangentially to the set of functions $h \in \mathbb{C}[0,1]$ 
with derivative $\phi^\prime: \mathbb{C}[0,1] \to \mathbb{R}$ given by $\phi^\prime(h) = h(x^*)$.
\end{lem}

We now show that under mild regularity conditions on $T$ (or, equivalently on $F$ and $G$) 
the estimator $\overline{m}_{T_n}$ is asymptotically normal.:
\begin{thm}
\label{thm:m-conv}
Let $F_n$ and $G_n$ be the empirical distribution functions of two independent random samples $X_1,\ldots,X_n$ and $Y_1,\ldots,Y_n$ from (absolutely continuous) distribution functions $F$ and $G$, respectively and let $T = G\circ F^{-}$. 
If $T$ is such that there exists a unique $p^*\in [0,1]$ with $T(p^*-) - p^* = \inf_{x\in[0,1]} T(x) - x$ and $F, G$ 
are such as in Corollary \ref{cor:FG-had-diff}, then for $T_n = G_n \circ F_n^{-}$,
$$	
  \sqrt{n} \big( \inf_x (T_n(x)-x) - \inf_x (T(x)-x) \big) 
$$
is asymptotically normal with mean $0$ and variance 
$$\frac{(G^\prime)^2(x_{p^*})}{(F^\prime)^2(x_{p^*})} \, p^*(1-p^*) + G(x_{p^*})(1-G(x_{p^*}))$$
where $x_{p^*} = F^{-}(p^*)$.
\end{thm}
\begin{proof} 
According to Donsker's Theorem (see \cite{vdV}) $(\mathbb{G}_{n,F}, \mathbb{G}_{n,G}) = \sqrt{n}(F_n-F, G_n-G)$ converges in distribution to $(\mathbb{G}_F, \mathbb{G}_G)$ in the space $\mathbb{D}[-\infty,\infty] \times \mathbb{D}[-\infty,\infty]$, for a pair of independent Brownian Bridges $\mathbb{G}_F$ and $\mathbb{G}_G$. The sample paths of the two limit processes are continuous, since both, $F$ and $G$, are continuous. By Corollary \ref{cor:FG-had-diff}, Lemma \ref{lem:m-had-diff} and the Chain Rule for Hadamard derivatives $\phi(F,G) = 1 + \inf_{x\in[0,1]} (G \circ F^{-}(x) - x)$ is Hadamard differentiable tangentially to the range of the limit processes. 
Applying the functional delta method yields that the sequence $\sqrt{n}(\inf_x (T_n(x)-x) - \inf_x (T(x)-x))$ is asymptotically equivalent to the derivative of $\phi$ evaluated at $(\mathbb{G}_{n,F}, \mathbb{G}_{n,G})$, i.e., to $-\frac{G^\prime(x_{p^*})}{F^\prime(x_{p^*})}\mathbb{G}_{n,F}(x_{p^*}) + \mathbb{G}_{n,G}(x_{p^*})$. 
Asymptotic normality now follows from the central limit theorem.
\end{proof}

Theorem \ref{thm:m-conv} considered uniqueness of the point attaining the infimum, the following final result 
considers the other extreme case where each point is a minimizer:
\begin{thm}
Let $F_n$ and $G_n$ be the empirical distribution functions of two independent random samples $X_1,\ldots,X_n$ and $Y_1,\ldots,Y_n$ and let $T = G\circ F^{-1}$. If $F, G$ are both $U(0,1)$ then $\sqrt{n}(\inf_x (T_n(x)-x) - \inf_x (T(x)-x))$ converges to $\inf_{t\in(0,1)} \sqrt{2} B_t$ (with $B_t$ being a standard Brownian Bridge) and thus has density $f(x) = -2x\exp(-x^2) \mathbf{1}_{(-\infty,0]}(x)$.
\end{thm}

\noindent The following final example illustrates Theorem \ref{thm:m-conv}. 

\begin{figure}[h!]
  \begin{center}
  \includegraphics[width = 12cm,angle=270]{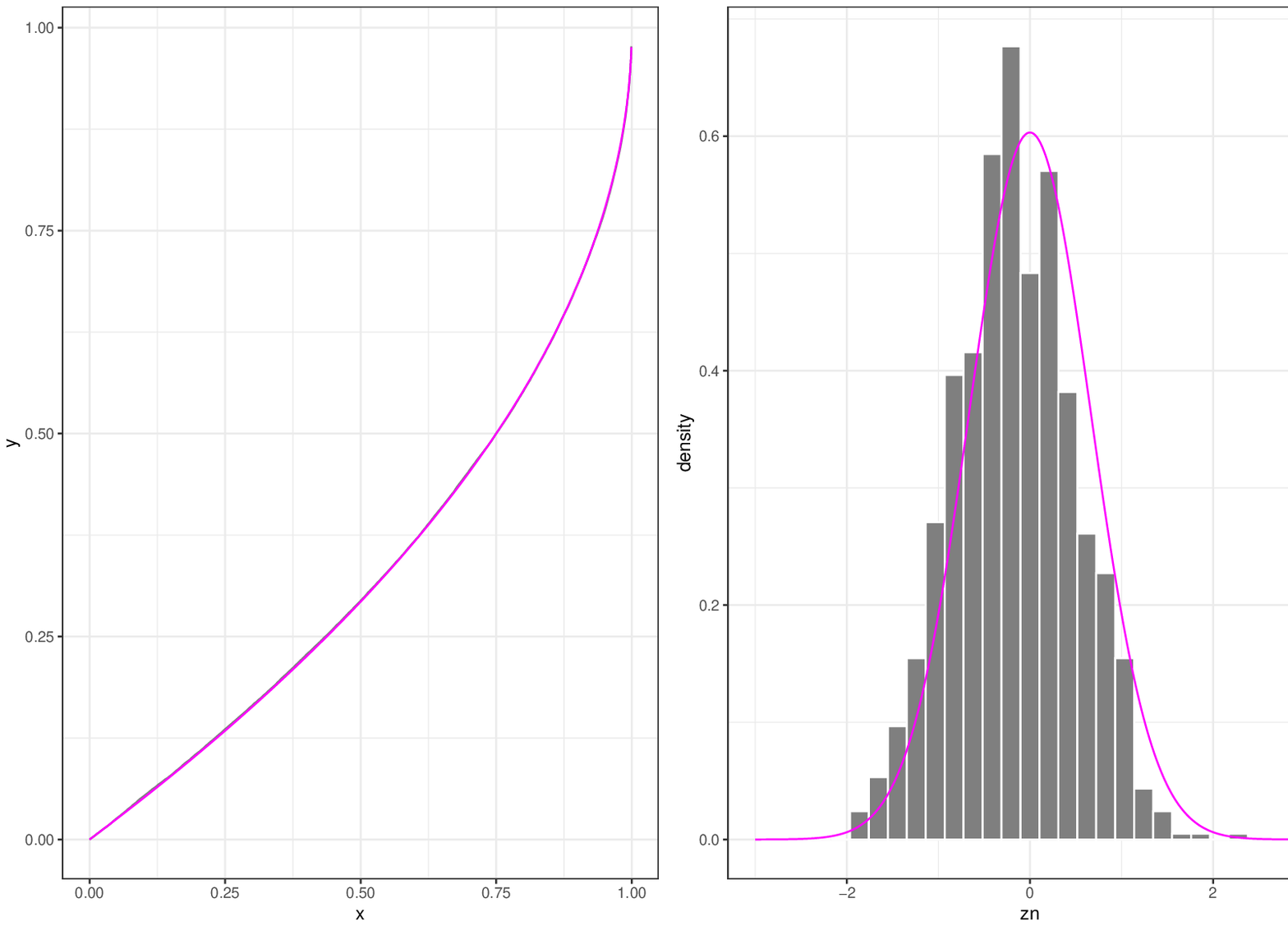}
        \caption{Left panel: $T$ (magenta) and $T_n$ with $n=100.000$ as considered in Example \ref{exponential_lifetimes};
        Right panel: Histogram of $R=1.000$ values of $Z_n$ with $n=100.000$ and density of $\mathcal{N}(0,\frac{7}{16})$.}\label{fig:asvar}
  \end{center}
\end{figure}

\begin{ex}[Example \ref{exponential_lifetimes} continued]
Consider the setting from Example \ref{exponential_lifetimes} for the case $\theta_1=2$ and $\theta_2=1$. 
Then it is straightforward to verify that all assumptions of Theorem \ref{thm:m-conv} are fulfilled, 
that $p^*=\frac{3}{4}$ is the unique minimizer, that $x_p=F^-(p^*)=\ln{2}$, and that
the asymptotic variance $\sigma^2$ is given by $\sigma^2=\frac{7}{16}$.
The right panel in Figure \ref{fig:asvar} depicts a histogram of $R=1.000$ samples of the random variable 
$Z_n:=\sqrt{n}(\overline{m}_{T_n}-\overline{m}_T)$ calculated by randomly drawing independent samples $X_1,\ldots,X_n$ 
and $Y_1,\ldots,Y_n$ from 
$X \sim Ex(\theta_1)$ and $Y \sim Ex(\theta_2)$ of size $n=100.000$, respectively.  
\end{ex}

\begin{rmk}
Based on simulations we conjecture that working with Bernstein approximations or splines it might be possible to derive 
strongly consistent estimators for $\overline{w}_T$ too. 
We plan to tackle this question in the near future.
\end{rmk}

\vspace{5mm}

\noindent \textbf{Acknowledgments}\\
The third and the fourth author gratefully acknowledge the support of the WISS 2025
project ‘IDA-lab Salzburg’ (20204-WISS/225/197-2019 and 20102-F1901166-KZP).



\begin{thebibliography}{00}
\bibitem{CCR} J. C\'arcamo, A. Cuevaz, L.-A. Rodr\'iguez: Directional differentiability for supremum-type functionals: statistical applications, 
   \emph{Bernoulli} \textbf{26}, 2143-2175 (2020). (see \url{https://arxiv.org/abs/1902.01136})

\bibitem{DFST} F. Durante, J. Fern\'andez S\'anchez, W. Trutschnig: On the singular components of a copula, 
   \emph{J. Appl. Probab.} \textbf{52}, 1175-1182 (2015).

\bibitem{DNO} W.F. Darsow, B. Nguyen, E.T. Olsen: Copulas and Markov processes, \emph{Illinois J. Math.} \textbf{36}, 600-642 (1992). 

\bibitem{DSS} F. Durante, P. Sarkoci, C. Sempi: Shuffles of copulas, \emph{J. Math. Anal. Appl.} \textbf{352}, 914-921 (2009).

\bibitem{DKQS} F. Durante, E.P. Klement, J. Quesada-Molina, P. Sarkoci: Remarks on Two Product-like Constructions for Copulas, \emph{Kybernetika} \textbf{43}, 235--244 (2007).

\bibitem{Du} F. Durante, C. Sempi: \emph{Principles of Copula Theory}, Chapman and Hall/CRC, 2015.

\bibitem{El} J. Elstrodt: \emph{Mass- und Integrationstheorie}, Springer, (1999).

\bibitem{EP} P. Embrechts, G. Puccetti: Bounds for functions of dependent risks, \emph{Finance and Stoch}  \textbf{10}, 341-352 (2006).

\bibitem{EP2} P. Embrechts, G. Puccetti: Bounds for functions of multivariate risks, \emph{J. Multivariate Anal.}  \textbf{97}, 526-547 (2006).

\bibitem{EH} P. Embrechts, M. Hofert: A note on generalized inverses, \emph{Math. Method. Oper. Res.} \textbf{77}, 423-432 (2013).

\bibitem{FST} J. Fern\a'andez S\a'anchez, W. Trutschnig: Conditioning based metrics on the space of multi\-variate copulas and 
         their interrelation with uniform and levelwise convergence and Iterated Function Systems,  
          \emph{Journal of Theoretical Probability} \textbf{28}, 1311-1336 (2015). 
          
\bibitem{FSTComm} J. Fern\a'andez S\a'anchez, W. Trutschnig: Some members of the class of (quasi-) copulas with given diagonal from the Markov kernel perspective, 
         Comm. Stat. A–Theor. \textbf{45}, 1508-1526 (2016).    
         
\bibitem{HS} E. Hewitt, K. Stromberg: \emph{Real and Abstract Analysis}, Springer Verlag, Berlin Heidelberg, (1965).

\bibitem{Ka} O. Kallenberg: \emph{Foundations of modern probability}, Springer Verlag, New York Berlin Heidelberg, (1997).

\bibitem{Kl} A. Klenke: \emph{Probability Theory - A Comprehensive Course}, Springer Verlag, Berlin Heidelberg, (2007).

\bibitem{Lan}  H.O.~Lancaster: Correlation and complete dependence of random variables, \emph{Ann. Math. Stat.} \textbf{34}, 1315-1321 (1963).
 
\bibitem{MS} J.F. Mai, M. Scherer: Simulating from the copula that generates the maximal probability for a joint default under given (inhomogeneous)
          marginals, \underline{in} \emph{Topics from the 7th International Workshop on Statistical Simulation} ed. V. Melas, S. Mignani, P. Monari, and L. Salmaso,
          Springer Proceedings in Mathematics \& Statistics \textbf{114}, pp. 333-341, 2014.    
         
\bibitem{MST} P. Mikus\'inski, H. Sherwood, M.D. Taylor: Shuffles of Min, \emph{Stochastica} \textbf{290}, 61-74 (1992).          
          
\bibitem{MFT} T. Mroz, S. Fuchs, W. Trutschnig: How simplifying and flexible is the simplifying assumption in 
         pair-copula constructions - analytic answers in dimension three and a glimpse beyond, 
          \emph{Electronic Journal of Statistics} \textbf{15}, 1951-1992 (2021).

\bibitem{Nel} R.B. Nelsen: An Introduction to Copulas, Springer, New York, (2006).          
    
\bibitem{Ru} W. Rudin: \emph{Real and Complex Analysis}, McGraw-Hill International Editions, Singapore, (1987).   

\bibitem{Rue} L. R\"uschendorf: Random Variables with Maximum Sums, Advances in Applied Probability, Vol. 14, No. 3, 623-632 (1982).

\bibitem{Ry} J.V. Ryff: Measure Preserving Transformations and Rearrangements, \emph{J. Math. Anal. Appl.}  \textbf{31}, 449-458 (1970). 

\bibitem{Th} H. Thorisson: \emph{Coupling, stationarity, and regeneration}, Probability and its Applications, Springer-Verlag, New York, (2000).

\bibitem{Tru} W. Trutschnig: On a strong metric on the space of copulas and its induced dependence measure, 
     \emph{J. Math. Anal. Appl.}  \textbf{384}, 690-705 (2011).

\bibitem{TFS} W. Trutschnig, J. Fern\'andez S\'anchez: Some results on shuffes of two-dimensional copulas, \emph{J. Stat. Plan. Infer.} \textbf{143}, 251-260 (2013).

\bibitem{vdV} A.W. van der Vaart: Asymptotic Statistics, Cambridge Series in Statistical and Probabilistic Mathematics, Cambridge University Press, (2000).

\bibitem{Zam} T. Zamfirescu: Most Monotone Functions are Singular, \emph{The American Mathematical Monthly} \textbf{88}(1), 
47-49 (1981)


\end{thebibliography}
\end{document}